\documentclass[12pt]{amsart}
\usepackage{amsmath,amssymb}

\newtheorem{theorem}{Theorem}[section]
\newtheorem{lemma}[theorem]{Lemma}
\newtheorem{keylemma}[theorem]{Key Lemma}
\newtheorem{proposition}[theorem]{Proposition}
\newtheorem{corollary}[theorem]{Corollary}

\theoremstyle{definition}
\newtheorem{definition}[theorem]{Definition}

\newtheorem{notation}[theorem]{Notation}

\renewcommand{\a}{\alpha}
\newcommand{\e}{\varepsilon}
\renewcommand{\d}{\delta}
\newcommand{\var}{\varphi}
\newcommand{\g}{\gamma}
\newcommand{\G}{\Gamma}
\renewcommand{\l}{\lambda}
\renewcommand{\o}{\omega}
\renewcommand{\O}{\Omega}
\newcommand{\s}{\sigma}
\renewcommand{\th}{\theta}
\newcommand{\Th}{\Theta}
\newcommand{\z}{\zeta}

\newcommand{\x}{\times}
\renewcommand{\i}{\infty}
\newcommand{\<}{\langle}
\renewcommand{\>}{\rangle}
\newcommand{\bw}{\bigwedge}
\newcommand{\p}{\partial}
\newcommand{\dc}{{\dot \cup}}

\newcommand{\bC}{{\mathbb C}}

\newcommand{\bR}{{\mathbb R}}
\newcommand{\bT}{{\mathbb T}}
\newcommand{\bZ}{{\mathbb Z}}

\newcommand{\cA}{{\mathcal A}}
\newcommand{\cB}{{\mathcal B}}
\newcommand{\cC}{{\mathcal C}}
\newcommand{\cD}{{\mathcal D}}
\newcommand{\cE}{{\mathcal E}}
\newcommand{\cH}{{\mathcal H}}
\newcommand{\cK}{{\mathcal K}}
\newcommand{\cL}{{\mathcal L}}
\newcommand{\cN}{{\mathcal N}}
\newcommand{\cO}{{\mathcal O}}
\newcommand{\cP}{{\mathcal P}}
\newcommand{\cT}{{\mathcal T}}
\newcommand{\cV}{{\mathcal V}}

\newcommand{\fG}{{\mathfrak g}}
\newcommand{\fH}{{\mathfrak h}}
\newcommand{\fM}{{\mathfrak m}}

\numberwithin{equation}{section}

\begin{document}

\title[Vector Bundles and Gromov--Hausdorff Distance]{Vector Bundles
and \\ Gromov--Hausdorff Distance }
\author{Marc A. Rieffel}
\address{Department of Mathematics\\
University of California\\
Berkeley, CA\ \ 94720-3840}
\email{rieffel@math.berkeley.edu}
\thanks{The research reported here was supported in part by National 
Science Foundation Grant DMS-0500501.}

\begin{abstract}
We show how to make precise the vague idea that for compact metric
spaces 
that are close together for Gromov--Hausdorff distance, suitable
vector 
bundles on one metric space will have counterpart vector bundles on
the 
other. Our approach employs the Lipschitz constants of
projection-valued functions that determine vector bundles.
We develop some computational techniques, and we illustrate our 
ideas with simple specific examples involving vector bundles on the
circle, 
the two-torus, the two-sphere, and finite metric spaces.  Our topic is
motivated by statements 
concerning ``monopole bundles'' over matrix algebras in the literature
of theoretical high-energy physics.
\end{abstract}

\maketitle

\section*{Introduction}

The purpose of this paper is to make precise the vague idea that if
two 
compact metric spaces are close together then there should be a
relationship 
between the vector bundles on the two spaces.  I was led to examine
this idea 
by statements in the theoretical high-energy physics literature 
\cite{Vlt1, BBV, GRS, Stn, BlI, Vlt2, CSW, GRS} stating, for example,
that 
for matrix algebras ``close'' to the two-sphere certain ``vector
bundles'' 
for the matrix algebras are {\em the} ``monopole bundles''
corresponding to 
the ordinary monopole bundles over the two-sphere.  I was able to make
sense 
of the statement that matrix algebras are close to the two-sphere
\cite{Rf73} 
by introducing \cite{Rf72} a definition of quantum metric spaces and
showing 
how to view the matrix algebras as such, and then by defining a
quantum 
Gromov--Hausdorff distance, which supplied a distance between the
matrix 
algebras and the two-sphere.  But I did not find in the literature any 
discussion  for ordinary metric spaces of a relationship between vector
bundles 
and ordinary Gromov--Hausdorff distance that I could then adapt to
the 
quantum setting.  The purpose of this paper is to develop such a
relationship 
for ordinary metric spaces.  I believe I have also found a path to such a 
relationship for quantum spaces, but it appears to be substantially
more 
complicated and indirect. (See \cite{Rf78}, where I have laid more
foundation for establishing that relationship, motivated by the
results of the present paper.) Thus I feel that it is worthwhile to
first explain 
separately this relationship just for ordinary metric spaces.  That is the
aim of the present paper.

In order to see the issues involved, let us consider a relatively
simple 
situation.  Consider, for example, the two-sphere $S^2$ of radius $1$
with 
its usual metric.  For some small $\e > 0$ let $F$ be a finite subset
of 
$S^2$ which is $\e$-dense in $S^2$, that is, every point of $S^2$ is
within 
distance $\e$ of a point of $F$.  (Matrix algebras are often
considered to 
be the algebras of ``functions'' on ``{\em quantum} finite sets''.)
Put on 
$F$ the metric from $S^2$.  Any vector bundle over $S^2$ restricts to
a 
vector bundle over $F$.  But all vector bundles over $F$ are trivial.
So, 
for example, the various inequivalent line bundles over $S^2$ will
restrict 
to equivalent line bundles over $F$.  At the topological level there
seems to 
be little more that one can say about this.  It is hard to see how one
can 
say that one line bundle over $F$ corresponds to the monopole bundle
on 
$S^2$, while a different one corresponds to the trivial line bundle on
$S^2$.

But vague intuition suggests that the restriction to $F$ of a
non-trivial 
line bundle over $S^2$ should somehow twist more rapidly than the
restriction 
to $F$ of a trivial bundle over $S^2$.  We want to make this intuition 
precise.  We need metrics in order to make sense of ``more rapidly''.
To see how 
to use metrics, let us drop back to the simpler example of the circle
and 
its simplest non-trivial vector bundle, the M\"obius-strip bundle.  If
you 
imagine explaining to a friend what this bundle is, using only your
hand, 
you will probably move your hand around in a circle, but twisting your
hand 
as it moves.  It makes sense to ask for the rate of twisting of your
hand 
with respect to arc-length along the circle.  To formulate this idea 
mathematically it seems necessary to view your hand as indicating a 
one-dimensional subspace twisting in the normal bundle to the circle 
within $\bR^3$.  This one-dimensional subspace can be specified by the 
orthogonal projection onto it.  More generally, it is well-known that
any 
vector bundle over a compact base space can be described up to
equivalence 
(in many ways) by means of a continuous projection-valued function on
the 
base space.  Given a metric on the base space, we can consider the
Lipschitz 
constant, $L(p)$, of a projection-valued function $p$.  We will see
that 
this provides a quite effective measure of how rapidly the bundle
determined 
by $p$ is twisting.  In fact, I have not seen any better way to
quantitatively discuss how rapidly a vector bundle twists.

We will find that for a given compact base metric space and a given
constant 
$K$ the set of projection-valued functions of given rank and size can
have 
several path components (even for a finite metric space).  These
components 
reflect different amounts and types of twisting.  We will find that
for 
examples such as an $\e$-dense $F \subset S^2$ and a $p$ defined on $F$,
and for 
more general Gromov--Hausdorff contexts, the constant $K$ can be
chosen in 
terms of $\e$ so that connected components of the set of $p$'s with 
$L(p) < K$ can only come from corresponding projection-valued
functions 
$q$ on $S^2$ which give isomorphic vector bundles.  That is, we have a 
uniqueness theorem (Theorems~\ref{th4.5} and \ref{th4.7}) in this
context.  
We also have an existence theorem, which states in quantitative terms
that 
if $\e L(p)$ is small enough, then there will exist a corresponding 
projection-valued function $q$ on the big space whose restriction is
$p$, 
such that $L(q)$ too is appropriately small.  We also have (and need) 
homotopy versions of these theorems (Theorems~\ref{th6.2} and
\ref{th6.5}). To give an indication of the nature of our results, we state the
following imprecise version of Theorem~\ref{th6.5}.

\begin{theorem}[Imprecise version of Theorem 6.5]
Let $X$ and $Y$ be compact metric spaces, and let $\rho$ be a metric on
their disjoint union that restricts to the given metrics on $X$ and $Y$, 
and for which the Hausdorff distance between them is less
than $\epsilon$. Let $p_0$ and $p_1$
be functions on $X$ whose values are projections in $M_n(\mathbb R)$,
so that they determine vector bundles on $X$.
There are positive constants $r_1, \ r_2, \ r_3$, depending
on $\epsilon$ and $n$,  such that the following holds. Assume that there is a
continuous path $p_t$ of projection-valued functions on $X$ going from $p_0$
to $p_1$ such that $L_X(p_t) < r_1$ for all $t$, where $L_X$ denotes Lipschitz
constant for the metric on $X$.  Then there exist functions $q_0$ and $q_1$
on $Y$ whose values are projections in $M_n(\mathbb R)$ such that
$L^\rho(p_j \oplus q_j) < r_2$ for $j = 0, 1$. Furthermore, for any such
$q_0$ and $q_1$ there is a continuous path $q_t$ of projection-valued
functions on $Y$ going from $q_0$ to $q_1$ such that $L_Y(q_t) < r_3$
for all $t$. In particular, the vector bundles determined by $q_0$ and $q_1$
are isomorphic. It is appropriate to view these vector bundles on $Y$ as corresponding
to the vector bundles determined by $p_0$ and $p_1$.  
\end{theorem}
The precise version of Theorem~\ref{th6.5} gives, among other things, 
formulas for $r_1, \ r_2, \ r_3$. 

When the metric spaces are sufficiently ``nice'',  standard facts in comparison 
geometry can be used to relate vector bundles on ones 
that are close together. For example, the contents
of section 3 of \cite{Pds} show that if $\rho$ is a metric on the
disjoint union $X \dc Y$ of compact spaces such that $X$ is
$\epsilon$-close to $Y$, and if $Y$ is suitably locally geometrically
$(n-1)$-connected where $n$ is the covering dimension of $X$,
then there is a map $f$ from $X$ to $Y$ such that $\rho ( x, f(x))$ 
is suitably small for all $x \in X$; furthermore any two such $f$'s are homotopic 
if $X$ is an absolute neighborhood retract. Such an $f$ can be
used to pull back vector bundles from $Y$ to $X$. However, the
techniques that we will use here do not require conditions of
finite dimension, local connectivity, and ANR on the compact
metric spaces that we consider. In particular, we very much want to permit
$Y$ to be finite, so not connected. (Full matrix
algebras are viewed as the algebras of functions on nonexistent ``quantum
finite sets''.) We also need to have good control of the 
Lipschitz constants of projection-valued functions for the vector bundles, and 
it is not clear to me to what extent techniques such as 
those in \cite{Pds} can be used to get quantitative estimates on
Lipschitz constants that are as strong as those that we obtain. Perhaps 
a combination of the techniques
could in favorable situations give stronger estimates than those obtained here.  
(It is interesting to 
recall that Vietoris first defined homology groups in the context of
compact 
metric spaces, using the metric in an essential way. See \cite{Rts}
and its 
references. His methods seem to be in the spirit of those used in the 
present paper, but I have not seen a specific relationship.) 

Let us mention some related ideas.  If one looks back at the
high-energy 
physics literature which discusses ``vector bundles'' on matrix
algebras 
that converge to bundles on the two-sphere or other spaces (see
references 
above), no definition of ``convergence'' is given, but what is noted
is that 
the formulas for the Chern classes of the ``vector bundles'' on the
matrix 
algebras appear to converge to the Chern classes of the limit bundle.  
However, given a convergent sequence of any kind, one can always
change 
any one given term of the sequence without affecting the convergence
of 
the sequence.  Thus it does not seem possible to use this approach to 
justify asserting that a given ``vector bundle'' on a given matrix
algebra 
is {\em the} monopole bundle.  As another approach, if one has a
sequence 
of ordinary compact metric spaces which converges to a given compact
metric 
space, and if one has vector bundles over all of these spaces, one can
try  
to put a compatible metric on the disjoint union of the metric spaces
so that the sequence converges to the limit for Hausdorff distance, in 
such a way that one can combine the vector bundles to form a
(continuous) vector bundle 
over this disjoint union.   If one can do this, then one can say that
the 
sequence of vector bundles converges to the bundle on the limit space.
This 
approach is taken in \cite{Hwk1, Hwk2}. But again it does not seem
possible 
to use this approach to justify saying that a given bundle on one
given 
space of the sequence is {\em the} bundle on that space which
corresponds to 
the bundle on the limit space.

The first two sections of this paper develop the basic properties of
the 
Lipschitz seminorms that we will need, while the third section relates
these 
seminorms to projections.  Section~\ref{sec4} contains our general
uniqueness 
theorems for extending vector bundles, while the next two sections
develop 
our existence theorems for extending vector bundles, and discuss their 
consequences.  The next nine sections examine simple examples
that illustrate 
our general theory. These examples involve vector bundles on the
circle, 
the two-torus and the two-sphere (including monopole bundles), as well
as on finite metric spaces.  We develop 
techniques for actually finding projections corresponding to given
vector 
bundles, and for calculating the Lipschitz constants of these
projections.  
In Section~\ref{sec10} we indicate how Chern classes can be used to
obtain 
lower bounds for the Lipschitz constants.  Sections \ref{sec14} and 
\ref{sec15} discuss vector bundles on compact homogeneous spaces in 
preparation for discussing monopole bundles in Section~\ref{sec16}.  
Our treatment is far from exhaustive, even for the simple spaces we 
consider.  The purpose of the examples we discuss is only to 
provide ``proof of concept''.  There is much more to be explored 
(and we say nothing here about non-compact spaces).

The audience I have had in mind when writing this paper consists of 
geometers and topologists.  For their convenience I have selectively 
included discussion of certain known facts from analysis that we need, 
rather than just giving references to the literature.

I developed an early part of the material reported here during a two-month
visit at the Institut Mittag-Leffler and a one-month visit at the
Institut des Hautes Etudes Scientifique in the Fall of 2003. I am very 
appreciative of the quite stimulating and enjoyable conditions
provided by 
both institutes.

\tableofcontents


\section{The setting}
\label{sec1}

Let $(X,\rho_X)$ and $(Y,\rho_Y)$ be two compact metric spaces.  
Suppose that the Gromov--Hausdorff distance \cite{Grm2, Bll, Ska}
between 
them is $< \e$.  Then, by definition, there exist a compact metric
space 
$(Z,\rho)$ and isometric inclusions of $X$ and $Y$ into $Z$ such that,
when 
$X$ and $Y$ are viewed as subsets of $Z$, their Hausdorff distance is 
$< \e$.  Our aim is to show that such an inclusion provides a
correspondence 
between certain vector bundles over $X$ and certain vector bundles
over $Y$.  
Since we can always cut $Z$ down to $X \cup Y$, the issues essentially
come 
down to relating vector bundles on $X$ to those on $Z$.  Thus most of
our 
technical development will take place in the setting of a given
compact 
metric space $(Z,\rho)$ and a closed subset $X$ of $Z$, with the
metric on 
$X$ being the restriction of $\rho$ to $X$.  To say that the Hausdorff 
distance from $X$ to $Z$ is $\leq \e$ is the same as saying that $X$
is 
$\e$-dense in $Z$, that is, every element of $Z$ is within distance
$\e$ of 
some point of $X$.  We will usually use this terminology
``$\e$-dense''.

We are interested in both real and complex vector bundles.  But we
will see 
that the bookkeeping for complex vector bundles is slightly more
complicated 
than that for real vector bundles.  For this reason our general
theoretical 
discussion, here and later, will be phrased primarily for complex
bundles, 
and we will usually only discuss real vector bundles in those places
where 
it is not entirely clear how our discussion for complex vector bundles
can 
be adapted to the case of real vector bundles.  We will let $C(X)$
denote 
the algebra of continuous complex-valued functions on a compact space
$X$, 
and we will let $C_{\bR}(X)$ denote the corresponding algebra of
real-valued 
functions.  We equip both algebras with the supremum norm.

We can only expect to deal with vector bundles up to the usual bundle 
isomorphism.  It is well known \cite{Aty, Blk, CMR, Krb, Rsb, WO} that, up to
isomorphism, 
every complex vector bundle on $X$ corresponds (in many ways) to a
continuous 
function, $p$, on $X$ whose values are projection operators.
Specifically, 
$p$ has values in an $n \x n$-matrix algebra $M_n(\bC)$ for some $n$;
for any 
$x \in X$ the fiber at $x$ of the vector bundle corresponding to $p$
will be 
the subspace $p_x(\bC^n)$ of $\bC^n$.  Each $p_x$ is idempotent 
$((p_x)^2 = p_x)$, and we can, and will, require that $p_x$ is 
self-adjoint $((p_x)^* = p_x)$.  (This corresponds to choosing a 
Hermitian metric on the vector bundle.)  We can equally well view $p$
as an 
element of the algebra $M_n(A)$ where $A = C(X)$.  Then we simply have 
$p^2 = p = p^*$. If $X$ is connected then the rank of $p$ must be
constant. In this case $p$ can be viewed as having values in a
Grassman manifold \cite{Hsm}, since one way to define a
Grassman manifold is as the space of all projections of a given
rank on a given vector space.

How then will the metric control vector bundles on $X$?  We will use 
$\rho$ to define a corresponding Lipschitz seminorm, $L = L^{\rho}$,
on 
$M_n(A)$, and then use $L(p)$ for control.  To define $L$ we must
first say 
that on $M_n(\bC)$ we use the usual operator norm obtained by viewing 
elements of $M_n(\bC)$ as operators on the inner-product space
$\bC^n$.  
We then view $M_n(A)$ as consisting of continuous functions on $X$
with 
values in $M_n(\bC)$, and define the norm on $M_n(A)$ by 
$\|a\| = \|a\|_{\i} = \sup\{\|a(x)\|: x \in X\}$.  In terms of the
norm 
on $M_n(\bC)$ we define $L = L^{\rho}$ on $M_n(A)$ by
\[
L^{\rho}(a) = \sup\{\|a(x)-a(y)\|/\rho(x,y): x,y \in X,\ x \ne y\}
\]
for $a \in M_n(A)$.  We note that $L$ can easily be discontinuous for
the 
operator norm, and we can easily have $L(a) = +\i$.  Of course,
viewing $A$ 
as $M_1(A)$, we have $L$ defined on $A$, and then the definition of
$L$ given 
above becomes the traditional definition of the Lipschitz seminorm.
Let 
$\cA = \{f \in A: L(f) < \i\}$, the $*$-subalgebra of ``Lipschitz 
functions''.  It is well-known \cite{Wvr2} and easily seen that $\cA$
is a 
dense $*$-subalgebra of $A$, where $*$ refers to complex conjugation.
It is 
also easily seen that for any $a \in M_n(A)$ we have $L(a) < \i$ if
and only 
if each entry of the matrix $a$ is in $\cA$, that is, $a \in
M_n(\cA)$.  Of 
course, $M_n(\cA)$ is a dense $*$-subalgebra of $M_n(A)$.  Also, one
should 
notice that for $a \in M_n(A)$ we have $L(a) = 0$ if and only if $a$
is a 
constant function with value in $M_n(\bC)$.  (Thus $L(a) = 0$ for more 
elements $a$ than just the scalar multiples of the identity element of 
$M_n(A)$.)  We will see later (Proposition~\ref{prop3.1}) that every
vector 
bundle is given by a $p$ in $M_n(\cA)$, not just in $M_n(A)$.  Then 
$L(p) < \i$.  We will control a vector bundle by the $L(p)$'s for the 
possible $p$'s that represent it.  In fact, we will see through
examples 
why it is reasonable that projections $p$ be the main focus of our
attention.

In the situation in which $X \subseteq Z$, let $B = C(Z)$, and let 
$\pi: B \to A$ be the surjective homomorphism consisting of
restricting 
functions on $Z$ to $X$.  We will let $\pi$ also denote the
corresponding 
$*$-homomorphism from $M_n(B)$ onto $M_n(A)$.  We will let $L$ denote
also 
the Lipschitz seminorms on $B$ and $M_n(B)$, though we may write $L_B$ 
(or $L_A$ for $A$) when this may be helpful.  We let $\cB$ denote the 
dense $*$-subalgebra of $B$ consisting of the functions on $Z$ for
which $L$ is finite.

If $q$ is a projection in $M_n(\cB)$, then $p = \pi(q)$ will be a
projection 
in $M_n(\cA)$, and the vector bundle associated to $p$ is easily seen
to be 
the restriction to $X$ of the bundle on $Z$ associated with $q$.  Thus
an 
important part of our task will be to show that under suitable
conditions 
in terms of $L(p)$ and the Hausdorff distance between $X$ and $Z$, we
can 
prove that given a projection $p \in M_n(\cA)$ there exists a
projection 
$q \in M_n(\cB)$ such that $\pi(q) = p$, and such that $L(q)$
satisfies 
estimates that imply that $q$ is unique up to homotopy.  (We will see 
that Lipschitz-controlled homotopy is the appropriate equivalence in our
setting.)  
We will also deal with the situation in which we have two projections
$p_1$ 
and $p_2$ in $M_n(\cB)$ which are homotopic.  (Homotopic projections
give 
isomorphic vector bundles --- see 1.7 and 4.2 of \cite{Blk}.)
For all of this we first need to show that $L$ has strong properties.


\section{Properties of the Lipschitz seminorm $L$}
\label{sec2}

For our present purposes it will be crucial that $L^{\rho}$ satisfies 
certain technical properties.

\begin{definition}
\label{def2.1}

Let $L$ be a seminorm, possibly discontinuous and possibly taking 
value $+\i$, on a normed unital algebra $A$ (for 
which $\|ab\| \le \|a\|\|b\|$ and $\|1\| = 1$).  
We will say that $L$ satisfies the {\em Leibniz} property 
(or ``is Leibniz'') if for every $a,b \in A$ we have 
\[
L(ab) \le L(a)\|b\| + \|a\|L(b).
\]
If, in addition, whenever $a^{-1}$ exists in $A$ we have
\[
L(a^{-1}) \le L(a)\|a^{-1}\|^2,
\]
then we will say that $L$ is {\em strongly Leibniz}.  
If $A$ has an involution and if $L(a^*) = L(a)$ for all $a \in A$, 
we say that $a$ is a {\em $*$-seminorm}.

We will say that $L$ is {\em lower-semicontinuous} if for one 
$r \in \bR$ with $r > 0$ (hence for every $r > 0$) the set 
$\{a \in A: L(a) \le r\}$ is a norm-closed subset of $A$.  
We say that $L$ is {\em closed} if for one $r \in \bR$ with $r > 0$ 
(hence for every $r > 0$) the set $\{a \in A: L(a) \le r\}$ is closed 
in the completion, ${\bar A}$, of $A$.  We will say that $L$ is 
{\em semi-finite} if $\{a \in A: L(a) < \i\}$ is dense in $A$.

\end{definition}

It is evident that if $L$ is closed then it is lower-semicontinuous.  
Every lower-semicontinuous $L$ extends to a closed seminorm 
(see proposition~$4.4$ of \cite{Rf66}).  Let $\cA = \{a \in A: L(a) <
\i\}$.  
We can always define a new norm, $\|\cdot\|_L$, on $\cA$ by
\[
\|a\|_L = \|a\| + L(a).
\]
If $L$ is closed, then $\cA$ is complete for $\|\cdot\|_L$, 
as can be seen by adapting the proof of proposition~$1.6.2$ of
\cite{Wvr2}.  
If $L$ is Leibniz, then $\|\cdot\|_L$ is a normed-algebra norm, 
i.e., $\|ab\|_L \le \|a\|_L\|b\|_L$.  Thus if $L$ is both closed and 
Leibniz, then $\cA$ is a Banach algebra for $\|\cdot\|_L$.  I do not 
know of an example of a finite Leibniz seminorm which is not strongly
Leibniz. (See the comments after definition 1.2 of \cite{Rf78}.)

\begin{proposition}
\label{prop2.2}

For a compact metric space $(X,\rho)$, let $A = C(X)$ or $C_{\bR}(X)$,
and 
let $L = L^\rho$ on $M_n(A)$ be defined as in the previous section.
Then 
$L$ on $M_n(A)$ is a closed semi-finite strongly-Leibniz $*$-seminorm.

\end{proposition}

\begin{proof}
That $L$ is  a $*$-seminorm is trivial to verify.  The semi-finiteness
was 
indicated in the previous section.  That $L$ on $A$ is Leibniz is
well-known 
and easily shown by the same device as is used to show that a first 
derivative has the related Leibniz property.  For $L$ on $M_n(A)$ the 
proof is the same except that one must keep terms in the proper order, 
respecting the non-commutativity of $M_n(A)$.  

That $L$ is strongly Leibniz seems not so well-known, so we include
the proof 
here.  Let $a \in M_n(A)$ and assume that $a$ is invertible.  If $L(a)
= \i$ 
the desired inequality is trivially true.  So assume that $L(a) < \i$.
For $x,y \in X$ with $x \ne y$ we have
\[
a(x)^{-1} - a(y)^{-1} = a(x)^{-1}(a(y)-a(x))a(y)^{-1},
\]
so that
\[
\|a(x)^{-1} - a(y)^{-1}\|/\rho(x,y) \le 
\|a^{-1}\|^2 \|a(x)-a(y)\|/\rho(x,y) \le \|a^{-1}\|^2L(a).
\]
Upon taking the supremum over $x$ and $y$ we obtain the desired
inequality.  
In particular, $a^{-1} \in M_n(\cA)$.

Finally, for any fixed $x$ and $y$ with $x \ne y$ clearly the function 
$a \mapsto \|a(x)-a(y)\|/\rho(x,y)$ is norm-continuous.  But 
$L$ is the supremum of these continuous functions as $x$ and $y$ range 
over $X$.  Thus $L$ is lower-semicontinuous on $M_n(A)$ (where $L$ may 
take value $+\i$).  Since $M_n(A)$ is complete for its norm, it
follows that $L$ is closed.
\end{proof}

There are some further properties of $L$ and $M_n(\cA)$ that will be 
essential to us in producing projections controlled by $L$.  For any 
unital algebra $A$ over $\bC$ and any $a \in A$, the spectrum, 
$\s(a)$, of $a$ is defined by
\[
\s(a) = \{\l \in \bC: (\l - a)^{-1} \mbox{ does not exist in } A\}
\]
(where we systematically write $\l$ instead of $\l 1_A$).  We make the 
analogous definition for algebras over $\bR$.  Note that for strongly 
Leibniz seminorms the ``strongly'' part implies that for $a \in \cA$ 
its spectrum in $\cA$ coincides with its spectrum in $A$.

For $f \in C(X)$ its spectrum coincides with the range of $f$.  Let 
$A = C(X)$ and let $a \in M_n(A)$.  Then $a^{-1}$ exists in $M_n(A)$ 
exactly if $a(x)^{-1}$ exists in $M_n(\bC)$ for all $x$.  From this
one 
sees that $\s(a) = \bigcup_{x \in X} \s(a(x))$.  Basic Banach-algebra 
arguments \cite{KdR, Rdn} show that $\s(a)$ is a closed bounded subset 
of $\bC$.

Let $f \in A = C(X)$, and let $\var$ be a holomorphic function defined
on an 
open neighborhood of $\s(f)$, that is, on the range of $f$.  Then the 
composition $\var \circ f$ is well-defined and in $A$.  We write it as 
$\var(f)$.  But we also need to define $\var(a)$ for $a \in M_n(A)$,
where 
now we ask that $\var$ be holomorphic in a neighborhood of $\s(a)$.
It is 
not immediately clear how to do this, but a basic Banach-algebra
technique 
(see 3.3 of \cite{KdR}, or \cite{Rdn}) does this by means of the
Cauchy 
integral formula, and this technique is called the
``holomorphic-function (or symbolic) 
calculus''.  In the standard way used for ordinary Cauchy integrals,
we let 
$\g$ be a collection of piecewise-smooth oriented closed curves in the
domain 
of $\var$ that surrounds $\s(a)$ but does not meet $\s(a)$, such that
$\var$ on $\s(a)$ is represented by its Cauchy integral using $\g$.
Then 
$z \mapsto (z-a)^{-1}$ will, on the range of $\g$, be a well-defined 
and continuous function with values in $M_n(A)$.  Thus we can define 
$\var(a)$ by
\[
\var(a) = \frac {1}{2\pi i} \int_{\g} \var(z)(z-a)^{-1}dz.
\]
For a fixed neighborhood of $\s(a)$ the mapping $\var \mapsto \var(a)$
is a 
unital homomorphism from the algebra of holomorphic functions on this 
neighborhood of $\s(a)$ into $M_n(A)$ \cite{KdR, Rdn}.

\begin{proposition}
\label{prop2.3}
With notation as above, if $a \in M_n(\cA)$ then $\var(a) \in
M_n(\cA)$.  
In fact,
\[
L(\var(a)) \le \left( \frac {1}{2\pi} 
\int_{\g} |\var(z)|d|z|\right)(M_{\g}(a))^2L(a)
\]
where $M_{\g}(a) = \max\{\|(z-a)^{-1}\|: z \in \mathrm{range}(\g)\}$.
\end{proposition}

\begin{proof}
Because $L$ is lower-semicontinuous, it can be brought within the
integral 
defining $\var(a)$, with the evident inequality.  (Think of
approximating 
the integral by Riemann sums.)  Because $L$ is strongly Leibniz, this
gives
\[
L(\var(a)) \le \frac {1}{2\pi} \int_{\g} |\var(z)| \ 
\|(z-a)^{-1}\|^2 \ L(a) \ d|z|.
\]
On using the definition of $M_{\g}(a)$ we obtain the desired
inequality.
\end{proof}

The commonly used terminology for the fact that if $a \in M_n(\cA)$
then 
$\var(a) \in M_n(\cA)$ is that ``$M_n(\cA)$ is closed under the 
holomorphic-function calculus of $M_n(A)$''.  See, e.g., section~$3.8$
of 
\cite{GVF}.  We remark that Schweitzer has shown \cite{Swtz} that if
$A$ is 
any unital $C^*$-algebra and if $\cA$ is a unital $*$-subalgebra which
is 
closed under the holomorphic-function calculus of $A$, then $M_n(\cA)$
is 
closed under the holomorphic-function calculus of $M_n(A)$.  Thus our 
Proposition 2.3 is a special case of Schweitzer's theorem, but it is 
good to see the above simple direct proof for our special case. 

The proof of the above proposition depends strongly on working over
$\bC$.  
But we will to some extent be able to apply it when working over
$\bR$, 
in the following way. (See also \cite{Gdl}.) Notice that
$M_n(C_{\bR}(X))$ 
is a $*$-subring of $M_n(C(X))$ which is closed under multiplication
by 
scalars in $\bR$.  We will accordingly speak of ``real
$C^*$-subrings'', 
etc.  For a real $*$-subring $A$ of a $C^*$-algebra $B$ we will say
that $A$ 
is closed under the holomorphic-function calculus for $B$ if for every 
$a \in A$ with $a = a^*$ (so that $\s_B(a) \subset \bR$) and for every 
$\bC$-valued function $\var$ holomorphic in a neighborhood of
$\s_B(a)$ 
and taking real values on $\s_B(a)$ we have $\var(a) \in A$, where 
$\var(a)$ is initially defined as before to be an element of $B$.

\begin{proposition}
\label{prop2.4}
Let $(X,\rho)$ be a compact metric space, let $A = C_{\bR}(X)$, let
$L$ 
and $\cA$ be as just above for $A$, and let $B = C(X)$.  Then
$M_n(\cA)$ 
is closed under the holomorphic-function calculus for $M_n(B)$, for
every 
$n$.  If $a \in M_n(\cA)$ with $a = a^*$, and if $\var$ is holomorphic
in 
an open neighborhood of $\s_B(a)$ and $\var(\s_B(a)) \subset \bR$,
then 
the estimate of Proposition~$\ref{prop2.3}$ for $L(\var(a))$ holds
here also.
\end{proposition}

\begin{proof}
Let $a \in M_n(\cA)$ with $a = a^*$, and let $\var$ be a $\bC$-valued 
function holomorphic in a neighborhood of $\s_B(a)$ with 
$\var(\s_B(a)) \subset \bR$.  The closed unital $*$-subalgebra over 
$\bC$ generated by \ $a$ \ in $M_n(C(X))$ is naturally isomorphic to 
$C(\s(a))$ by one version of the spectral theorem, i.e. 
the ``continuous-function calculus'' (see 1.2.4 of \cite{KdR}, or page
8 
of \cite{Dvd}, or VI.6.3 of \cite{FlD}), with \ $a$ \ represented by
the 
function ${\hat a}$ given by ${\hat a}(r) = r$ for $r \in \s(a)$.  
(Here $\s(a)$ is the spectrum of \ $a$ \ as an element of
$M_n(C(X))$.)  
Then for any continuous $\bC$-valued function $\psi$ on $\s(a)$ we can 
form $\psi \circ {\hat a}$, and this will correspond to an element of 
$M_n(C(X))$, which we denote by $\psi(a)$.  If $\psi(\s(a)) \subset
\bR$, 
then by approximating $\psi$ on $\s(a)$ uniformly by polynomials over
$\bR$ 
it is clear that $\psi(a) \in M_n(C_{\bR}(X))$.  But if $\th$ is a
polynomial 
over $\bC$, then one can check that $\th(a)$ obtained by substituting
$a$ 
into $\th$ coincides with defining $\th(a)$ via the
holomorphic-function 
calculus, and that if $\th$ is over $\bR$ then $\th(a) \in M_n(A)$
either 
way.  From this one can check that $\var(a)$ defined via 
$\var \circ {\hat a}$ coincides with its definition via the 
holomorphic-function calculus, and $\var(a) \in M_n(C_{\bR}(X))$ 
since $\var(\s(a)) \subset \bR$.  But from Proposition~\ref{prop2.3} 
we see that $L(\var(a)) < \i$ so that $\var(a) \in M_n(\cA)$.  
The inequality for $L(\var(a))$ then follows from
Proposition~\ref{prop2.3}.
\end{proof}

Most of the examples that we discuss later involve manifolds, and for 
manifolds it is useful to be able to apply calculus to our
considerations.  
One tool relating calculus to our Lipschitz context, which we will 
find useful later, is given by the following proposition.

\begin{proposition}
\label{prop2.5}
Let $A$ be a unital Banach algebra and suppose that $L$ is a 
lower-semicontinuous seminorm on $A$.  Let $\a$ be a (strongly
continuous) 
action of a Lie group $G$ on $A$, and assume that this action is by 
isometries for $L$, that is, $L(\a_x(a)) = L(a)$ for all $a \in A$ 
and $x \in G$.  Let $A^{\i}$ denote the algebra of smooth elements of 
$A$ for the action $\a$. (It is closed under the holomorphic-function 
calculus of $A$ --- see proposition 3.45 of  \cite{GVF}).  Then for
any 
$a \in A$ and any $\e > 0$ there is a \ $b \in A^{\i}$ such that 
$\|b\| \le \|a\|$, $L(b) \le L(a)$, and $\|b-a\| < \e$.
\end{proposition}

\begin{proof}
We use the usual smoothing argument (see, e.g., section $0.2$ 
of \cite{Tyl}).  Let $f \in C_c^{\i}(G)$ and assume further that 
$f \ge 0$ and $\int_G f(x)dx = 1$ for left Haar measure on $G$.  
For a given $a \in A$ set
\[
b = \int_G f(x)\a_x(a)dx.
\]
Standard simple arguments show that $b \in A^{\i}$, that $\|b\| \le
\|a\|$, 
and that $\|b-a\| < \e$ if $f$ is supported sufficiently closely to
the 
identity element of $G$.  But by the lower semicontinuity of $L$ we
have
\begin{eqnarray*}
L(b) &= &L\left( \int_G f(x)\a_x(a)dx\right) \\
&\le &\int_G f(x)L(\a_x(a))dx = \int_G f(x)L(a)dx = L(a).
\end{eqnarray*}
\end{proof}


\section{Projections and Lipschitz seminorms}
\label{sec3}

To control a vector bundle we will use bounds on $L(p)$ for
projections 
representing the vector bundle, since we take $L(p)$ as a measure of
how 
rapidly a vector bundle twists.  For this to work well we first need
to know that we can always find representing projections $p$ such that
$L(p) < \i$.  
After showing this, we will establish a related fact for homotopies
between projections. A well-known fact that is important for all we do
here is that if two projections are homotopic then the corresponding
vector bundles are isomorphic. See, for example, section 6 of chapter
1 of \cite{Rsb}, especially corollary 1.6.12.

The fact that we can find $p$'s with $L(p) < \i$ is a special case of
a general well-known \cite{Krb} fact about any dense $*$-subalgebra
closed 
under the holomorphic-function calculus in any unital $C^*$-algebra 
(see, e.g., section~$3.8$ of \cite{GVF}).  We sketch the proof here in 
several steps, both for the reader's convenience and because we will
need similar arguments later. The normed $*$-algebras $M_n(C(Z))$ are
examples of $C^*$-algebras, and whenever we write ``$C^*$-algebra''
readers can just think of $M_n(C(Z))$ if they prefer.

\begin{proposition}
\label{prop3.1}
Let $A$ be a unital $C^*$-algebra, and let $\cA$ be a dense
$*$-subalgebra 
closed under the holomorphic-function calculus in $A$.  Let $p$ be a 
projection in $A$.  Then for any $\d > 0$ there is a projection $p_1$
in 
$\cA$ such that $\|p-p_1\| < \d$.  If $\d < 1$ then $p_1$ is homotopic
to 
$p$ through projections in $A$.  The same conclusions hold if $A$ is a
real 
$C^*$-subring of a $C^*$-algebra, as in the setting of 
Proposition~$\ref{prop2.4}$.
\end{proposition}

\begin{proof}
For the first part we can assume that $\d < 1/2$.  Since $\cA$ is
dense in 
$A$ and $p = p^*$ we can find $a \in \cA$ such that $a^* = a$ and 
$\|p-a\| < \d < 1/2$.  We will use the next step several times again, 
so we state it as:

\begin{lemma}
\label{lem3.2}
If $p$ is a projection in $A$ and if $a \in A$ with 
$a = a^*$ and $\|p-a\| < \d$, then
\[
\s(a) \subset [-\d,\d] \cup [1-\d,1+\d].
\]
\end{lemma}

\begin{proof}
(e.g., lemma~$2.2.3$ of \cite{RLL}) Note that 
$\s(p) \subseteq \{0, 1\}$. Let $\l \in \bR$ with $\l \ne 0,1$ so that 
$\l-p$ is invertible.  Then
\[
1 \ - \ (\l - p)^{-1}(\l - a) \ = \ (\l - p)^{-1}(a - p) .
\]
Thus if $\|(\l - p)^{-1}(a - p)\| < 1$ then $(\l - p)^{-1}(\l - a)$ is
invertible 
(by Neumann series, i.e., geometric series, converging in $A$), and so
$\l - a$ is left invertible.  But
\[
\|(\l-p)^{-1}(a \ - \ p)\| \le \|(\l-p)^{-1}\| \|a \ - \ p\| < 
\|(\l-p)^{-1}\|\d,
\]
and $\|(\l-p)^{-1}\| = \max\{|\l|^{-1},|1-\l|^{-1}\}$.  
Thus if $\d <|\l|$ 
and $\d < |1-\l|$ then $\l -a$ is left invertible. A small variation of this argument shows that $\l -a$ is also right invertible. Thus 
$\l \notin \s(a)$.  This yields the desired
result, since $\s(a) \subset \bR$. 
\end{proof}

We continue sketching the proof of Proposition~\ref{prop3.1}.  Let
$\chi$ be 
defined on $\bC$ by $\chi(z) = 0$ if $\mathrm{Re}(z) \le 1/2$, while 
$\chi(z) = 1$ if $\mathrm{Re}(z) > 1/2$.  Since $\d < 1/2 < 1 - \d$, 
we see from Lemma~\ref{lem3.2} that $\chi$ is holomorphic in an open 
neighborhood of $\s(a)$, and thus $\chi(a)$ is defined and is in
$\cA$.  
Since $\var \mapsto \var(a)$ is a homomorphism, $p_1 = \chi(a)$ is a 
projection in $\cA$.  By considering the continuous-function calculus 
\cite{KdR, Dvd, FlD} it is easily seen that $\|p_1 - a\| < \d$, and so 
$\|p-p_1\| < 2\d$.  Replacing $\d$ by $\d/2$ everywhere above, we
obtain the desired result.

The proof that $p$ and $p_1$ are homotopic is a simpler version of the
proof of our next proposition, below.

On looking at Proposition~\ref{prop2.4} and its proof one can see
easily 
how to adapt the above discussion to treat a real $C^*$-subring of a 
$C^*$-algebra.
\end{proof}

When we apply Proposition~\ref{prop3.1} to $M_n(C(X))$ or
$M_n(C_{\bR}(X))$ 
we can interpret it as saying that every vector bundle over $X$ has a 
Lipschitz structure with respect to $\rho$.

\begin{proposition}
\label{prop3.3}
Let $A$ and $\cA$ be as in Proposition~\ref{prop3.1}, and let $L$ be a 
lower-semicontinuous strongly-Leibniz $*$-seminorm on $\cA$.  Let
$p_0$ and 
$p_1$ be two projections in $A$.  
Suppose that $\|p_0-p_1\| \le \d < 1$, so 
that there is a norm-continuous path, $t \mapsto p_t$, of projections
in $A$ 
going from $p_0$ to $p_1$ \cite{Blk, RLL}.  If $p_0$ and $p_1$ are in
$\cA$ 
then we can arrange that $p_t \in \cA$ and that
\[
L(p_t) \le (1-\d)^{-1} \max\{L(p_0),L(p_1)\}
\]
for every $t$.  The same conclusions hold if $A$ is a real
$C^*$-subring of a $C^*$-algebra. 
\end{proposition}

\begin{proof}
We begin the proof in a standard way (e.g. proposition 2.2.4
of\cite{RLL}) 
by setting $a_t = (1-t)p_0 + tp_1$ for $t \in [0,1]$.  
If $0 \le t \le 1/2$ 
then $\|a_t-p_0\| \le \d/2$, while if $1/2 \le t \le 1$ then 
$\|a_t-p_1\| \le \d/2$.  From Lemma~\ref{lem3.2} and the evident facts 
that $a_t$ is positive and $\|a_t\| \le 1$, it follows that 
$\s(a_t) \subseteq [0,\d/2] \cup [1-\d/2,1]$.  Much as in the proof of 
Proposition \ref{prop3.1}, define $\chi$ on $\bC$ 
by setting $\chi(z) = 0$ 
if $\mathrm{Re}(z) \le (1+\d)/4$, and $\chi(z) = 1$ if 
$\mathrm{Re}(z) > (1+\d)/4$.  Since $\d < 1$ we see that $\chi$ is 
holomorphic in an open neighborhood of $\s(a_t)$ for each $t$.  

Instead of using the curve around $[1 - \d/2, 1]$ that we used in our earlier versions of this paper, 
we use the family, $\{\g_s\}$, of curves used by Hanfeng Li in his proof 
of proposition 3.1 of \cite{Li}. This gives substantially improved estimates compared
to those in our earlier versions. For the reader's convenience we include the details from his
proof here. For $s > 0$ we let $\g_s$ be the oriented curve that traces counter-clockwise
at unit speed the boundary of the rectangle with vertices $(1/2)-si, \ 1 + s - si, \
 1 + s + si, \ (1/2) + si$.
Notice that because of where 
we have chosen the line of discontinuity for $\chi$, the curves $\g_s$
lie in 
the domain where $\chi$ is holomorphic.  We would also usually choose
a curve 
around $[0,\d/2]$, but since $\chi = 0$ near there, this is
unnecessary. For any given $t$ we now set
\[
p_t = \chi(a_t) = \frac {1}{2\pi i} \int_{\g_s} \chi(z)(z-a_t)^{-1}dz.
\]
This integral is independent of $s$ for the usual reasons for
holomorphic functions.
Then as in the proof of Proposition~\ref{prop3.1} we see that $p_t$ 
is a projection, and that each $p_t$ is in $\cA$ if $p_0$ and $p_1$
are.  
It is not difficult to see then that $t \mapsto p_t$ is a continuous
path from $p_0$ to $p_1$.

We now estimate $L(p_t)$.   For the same reasons 
as given in the proof of Proposition~\ref{prop2.3} we have
\[
L(p_t) \le \frac {L(a_t)}{2\pi} 
\left( \int_{\g_s} d|z|    \ \| (z - a_t)^{-1}\|^2\right) \  .
\]
Let $\g^1_s$ be the segment of $\g_s$ given by $\g^1_s(r) = (1/2) + r - si$ for $0 \leq r \leq (1/2) + s$. For $z$ on $\g^1_s$, because $a_t$ is self-adjoint, we have $\|(z - a_t)^{-1}\| \leq s^{-1}$. Thus
\[
 \int_{\g^1_s} d|z|    \ \| (z - a_t)^{-1}\|^2 \ \leq \ \int_0^{(1/2)+s} s^{-2} \ dr \ = \ s^{-2}((1/2) + s) \ .
\]
This same estimate applies to the segment $\g^3_s$ of $\g_s$ going from $(1/2) +s +si$
to $(1/2) + si$. Notice that these estimates show that as $s$ goes to $+\infty$ these integrals go to 0.
Now let $\g^2_s$ be the segment of $\g_s$ given by $\g^2_s(r) = 1 + s + ri$ for $-s \leq r \leq s$. 
For $z$ on $\g^2_s$ we again have the estimate $\|(z-a_t)^{-1}\| \leq s^{-1}$. Thus
\[
 \int_{\g^2_s} d|z|    \ \| (z - a_t)^{-1}\|^2 \ \leq \ \int_{-s}^{s} s^{-2} \ dr \ = 2s^{-1} \ .
\] 
Notice again that this integral goes to 0 as $s$ goes to $+\infty$. Finally, let $\g^4_s$
be the segment of $\g_s$ given by $\g^4_s(r) = (1/2) - ri$ for $-s \leq r \leq s$.
For $z$ on $\g^4_s$ we have the estimate $\|(z-a_t)^{-1}\|^2 \leq ((1/2) - \d/2)^2 + r^2)^{-1}$.
Thus
\[
\begin{aligned}
\int_{\g^4_s} d|z|    \ \| (z - a_t)^{-1}\|^2 \ &\leq \ \int_{-\infty}^\infty  ((1/2) - \d/2)^2 + r^2)^{-1} \ dr  \\
& = \ \pi((1/2) - \d/2)^{-1} = 2 \pi(1-\d)^{-1}.
\end{aligned}
\]
Taking the sum over the 4 segments and letting $s$ go to $+\infty$, we obtain
\[
L(p_t) \le L(a_t) (1-\d)^{-1} .
\]
Since clearly $L(a_t) \leq \max\{L(p_0),L(p_1)\}$, we obtain the desired estimate.

It is easy to see how to adapt the above argument to the case of a
real $C^*$-subring of a $C^*$-algebra.
\end{proof}

By using the same techniques, or by applying directly proposition 3.1 of \cite{Li},
we can obtain the following continuation
of Proposition~\ref{prop2.5}:

\begin{proposition}
\label{prop3.4}
Let $A$ be a unital $C^*$-algebra and let $L$ be a closed
strongly-Leibniz 
$*$-seminorm on $A$.  Let $\cA$ be its subalgebra of Lipschitz
elements.  
Let $\a$ be an action of a Lie group $G$ on $A$ by isometries for $L$, 
and let $A^{\i}$ be the subalgebra of smooth elements for $\a$.  Then 
for any projection $p$ in $\cA$ and any 
$\d > 0$ there is a projection $p_1$ in $A^{\i}$ such that 
$\|p-p_1\| < \d$ and $L(p_1) \le (1-2\d)^{-1}L(p)$.  In particular, 
if $\d < 1/2$ then $p$ and $p_1$ are homotopic through projections in
$\cA$.  
The same conclusion holds if $A$ is a real $C^*$-subring of a
$C^*$-algebra.
\end{proposition}

\begin{proof}
According to Proposition~\ref{prop2.5}, given a projection $p \in \cA$
and  
a \ $\d > 0$ we can find a \ $b \in A^{\i}$ such that  
$\|b\| \leq \|p\| = 1$, $L(b) \le L(p)$, and $\|p-b\| < \d$.  
In fact, examination of the proof of Proposition~\ref{prop2.5} 
shows that we can also assume that $b^*=b$, and even more 
that $b \ge 0$, since $p \ge 0$.  
Thus $\s(b) \subseteq [0,1]$.  Then from Lemma~\ref{lem3.2} we
conclude 
that $\s(b) \subseteq [0,\d] \cup [1-\d,1]$.  We can assume that $\d < 1/2$, 
so that these intervals are disjoint. Then from proposition 3.1 of
\cite{Li} we obtain $L(p_1) \le (1-2\d)^{-1}L(p)$, as desired.
\end{proof}

This proposition is related to the main theorem of \cite{Li}, which
answers a question that I asked in an earlier version of this paper. We
state Li's theorem here, since we will use it later.

\begin{theorem}
\label{th3.5}
Let $M$ be a closed connected compact Riemannian manifold, 
equipped with its usual metric from its Riemannian metric, and let
$A$ be a real C*-subring of a C*-algebra. For any projection $p \in C(M, A)$
and any $\epsilon > 0$ there exists a projection $q \in C^\infty(M,A)$ such
that $\|p-q\|_\infty < \epsilon$
and $L(q) \leq L(p) + \epsilon$
\end{theorem}


\section{The uniqueness of extended vector bundles}
\label{sec4}

Let $(Z,\rho)$ be a compact metric space, and let $X$ be a closed
non-empty 
subset of $Z$.  As before, we equip $X$ with the restriction to it of
$\rho$, 
and we let $A = C(X)$ and $B = C(Z)$.  Let $L = L^{\rho}$ be defined
as 
earlier on $A$ and $B$, and also on $M_n(A)$ and $M_n(B)$, and let
$\cA$ 
and $\cB$ denote the dense subalgebras of Lipschitz functions.  
We let $\pi: M_n(B) \to M_n(A)$ denote the surjective $*$-homomorphism 
consisting of restricting functions from $Z$ to $X$.  It is easily
seen 
that for $b \in M_n(B)$ we have $L(\pi(b)) \le L(b)$, so that 
$\pi(M_n(\cB)) \subseteq M_n(\cA)$.  We say that a projection 
$q \in M_n(B)$ extends a projection $p \in M_n(A)$ if $\pi(q) = p$.  
This corresponds exactly to extension for the corresponding vector
bundles.

As mentioned earlier, if two projections, $q_1$ and $q_2$, 
in $M_n(B)$ are homotopic through projections in $M_n(B)$, then the 
corresponding vector bundles are isomorphic.  Thus we seek conditions 
such that if $p$ is a projection in $M_n(A)$ and if $q_0$ and $q_1$
are 
projections in $M_n(B)$ such that $\pi(q_0) = p = \pi(q_1)$, then
$q_0$ 
and $q_1$ are homotopic.  Simple examples show that this need not hold 
without further conditions.  As indicated in Section~\ref{sec1}, our 
results will depend on the Hausdorff distance between $X$ and $Z$, 
and we will express this by supposing that $X$ is $\e$-dense in $Z$.  

\begin{keylemma}
\label{keylem4.1}
Suppose that $X$ is $\e$-dense in $Z$.  Then for any $b \in M_n(B)$ we
have
\[
\|b\| \le \|\pi(b)\| + \e L(b).
\]
\end{keylemma}

\begin{proof}
Let $z \in Z$ be given.  Then there is an $x \in X$ such that 
$\rho(z,x) \le \e$.  Thus
\[
\|b(z)\| \le \|b(x)\| + \|b(z)-b(x)\| \le \|\pi(b)\| + \e L(b).
\]
\end{proof}

When we use the word ``path'' in discussing homotopies we will almost 
always mean a continuous function whose domain is the interval
$[0,1]$.

\begin{theorem}
\label{th4.2}
Suppose that $X$ is $\e$-dense in $Z$.  Let $p_0$ and $p_1$ be
projections 
in $M_n(\cA)$, and let $q_0$ and $q_1$ be projections in $M_n(\cB)$
such 
that $\pi(q_0) = p_0$ and $\pi(q_1) = p_1$.  Set
\[
\d = \|p_0-p_1\| + \e(L(q_0) + L(q_1)).
\]
If $\d < 1$, then there is a path, $t \mapsto q_t$, through
projections 
in $M_n(\cB)$, from $q_0$ to $q_1$, such that
\[
L(q_t) \le (1-\d)^{-1} \max\{L(q_0),L(q_1)\}
\]
for all $t \in [0,1]$.  The same conclusion holds if $A = C_{\bR}(X)$,
etc.
\end{theorem}

\begin{proof}
From Key Lemma~\ref{keylem4.1} we see that
\begin{eqnarray*}
\|q_0-q_1\| &\le &\|\pi(q_0-q_1)\| + \e L(q_0-q_1) \\
&\le &\|p_0-p_1\| + \e(L(q_0) + L(q_1)) = \d.
\end{eqnarray*}
Assume now that $\d < 1$.  Then according to Proposition~\ref{prop3.3} 
applied to $q_0$ and $q_1$ for $\cA = \cB$, there is a path $t \to q_t$ from $q_0$ to 
$q_1$ with the stated properties.
\end{proof}

Under the conditions of the above theorem, $t \mapsto \pi(q_t)$ will be
a continuous path from $p_0$ to $p_1$ through projections in
$M_n(\cA)$.  
Thus the vector bundles corresponding to $p_0$ and $p_1$ will be
isomorphic, 
as will their lifts corresponding to $q_0$ and $q_1$.

Notice that the bound on $L(q_t)$ stated in the above theorem is 
independent of $n$.  This is in contrast to the existence results 
that we will obtain in Section~\ref{sec6}. If $p_0 = p_1$ above then
we can obtain some additional information:

\begin{proposition}
\label{prop4.3}
Let $p \in M_n(\cA)$, and let $q_0$ and $q_1$ be projections in
$M_n(\cB)$ 
such that $\pi(q_0) = p = \pi(q_1)$.  If $\e L(q_0) < 1/2$ and 
$\e L(q_1) < 1/2$, then there is a path, $t \to q_t$, through
projections 
in $M_n(\cB)$, from $q_0$ to $q_1$, such that $\pi(q_t) = p$ and 
\[
L(q_t) \le (1-\d)^{-1} \max\{L(q_0),L(q_1)\}
\]
for all $t$, where $\d = \e(L(q_0) + L(q_1))$.  The same conclusion
holds 
when $A = C_{\bR}(X)$, etc.
\end{proposition}

\begin{proof}
The proof follows the same lines as the proof of Theorem~\ref{th4.2}, 
except that now we must show that $\pi(q_t) = p$ for all $t$.  
For this, much as in the proof of 
Proposition~\ref{prop3.3} we set $a_t = (1-t)q_0 + tq_1$ .  It is
clear that $\pi(a_t) = p$ for all $t$.  
From the proof of Proposition~\ref{prop3.3} we have 
$q_t = \frac {1}{2\pi i} \int_{\g} (z-a_t)^{-1}dz$, 
and it follows easily that $\pi(q_t) = p$ for all $t$.
\end{proof}

We can combine the above results to obtain some information that does
not depend on $p_0$ and $p_1$ being close together:

\begin{corollary}
\label{cor4.4}
Let $p_0$ and $p_1$ be projections in $M_n(\cA)$, and let $q_0$ and 
$q_1$ be projections in $M_n(\cB)$ such that $\pi(q_0) = p_0$ and 
$\pi(q_1) = p_1$.  Let $N$ be a constant such that $L(q_j) \le N$ for 
$j = 0,1$.  Assume further that there is a path $p$ from $p_0$ to 
$p_1$ such that for each $t$ there is a projection ${\tilde q}_t$ in 
$M_n(\cB)$ such that $\pi({\tilde q}_t) = p_t$ 
and $L({\tilde q}_t) \le N$.  
Assume that $2 \e N < 1$, and pick $\d$ such that $2 \e N < \d < 1$.   
Then there is a continuous path $t \mapsto q_t$ of projections in 
$M_n(\cB)$ going from $q_0$ to $q_1$ such that
\[
L(q_t) \le (1-\d)^{-1}N
\]
for each $t$.  (But we may not have $\pi(q_t) = p_t$ for all $t$.)  
The same conclusion holds when $A = C_{\bR}(X)$, etc.
\end{corollary}

\begin{proof}
Pick a finite increasing sequence $\{t_i\}_{i=0}^k$, of points in
$[0,1]$ such that $t_0 = 0$, $t_k = 1$, 
and $\|p_{t_{i+1}}-p_{t_i}\| \le \d - 2\e N$ 
for each $i$ for $0 \le i \le k - 1$.  For each such $i$ pick a
projection $q'_i$ in $M_n(\cB)$ such that $\pi(q'_i) = p_{t_i}$ and
$L(q'_i) \le N$, 
with $q'_0 = q_0$ and $q'_k = q_1$.  Then, according to Key Lemma 
\ref{keylem4.1}, for each such $i$ we have 
\[
\|q'_{i+1} - q'_i\| \le \|p_{t_{i+1}} - p_{t_i}\| + 
\e(L(q'_{i+1})+L(q'_i)) \le \d - 2 \e N + 2 \e N = \d.
\]
According to Theorem~\ref{th4.2}, for each $i$ there is a continuous
path 
of projections in $M_n(\cB)$, going from $q'_i$ to $q'_{i+1}$, all of
whose 
Lipschitz norms are no greater than $(1-\d)^{-1}N$.  We can then 
concatenate these paths in the usual way to obtain the desired path 
$t \to q_t$.  (Notice that the function $t \mapsto L(q_t)$ need not be 
continuous, and we need not have $\pi(q_t) = p_t$.)
\end{proof}

Let us now see what consequences the above uniqueness results have for 
metric spaces that are close together.  In doing this it seems
simplest to 
notice that in defining Gromov--Hausdorff distance between two compact 
metric spaces it is sufficient to consider their {\em disjoint}
isometric 
embeddings into other metric spaces.  Since we can always then cut
down to 
their union, it suffices to take $Z = X  \dc  Y$, where this denotes
the disjoint 
union.  Thus it suffices to consider metrics $\rho$ on $X  \dc  Y$
whose 
restrictions to $X$ and $Y$ are their given metrics $\rho_X$ and
$\rho_Y$.  
We will write $\mathrm{dist}_H^{\rho}(X,Y) < \e$ to signify that the 
Hausdorff distance between $X$ and $Y$ in $Z = X  \dc  Y$ for $\rho$
is less than $\e$.

From $Z = X  \dc  Y$ we have $C(Z) = C(X) \oplus C(Y)$ as $*$-Banach 
algebras.  As before let $A = C(X)$, $B = C(Z)$, etc., and now 
also let $D = C(Y)$, with subalgebra of Lipschitz elements $\cD$.  
A projection in $M_n(B)$ will now be of the form $p \oplus q$ where 
$p$ and $q$ are projections in $M_n(A)$ and $M_n(D)$ respectively.  
Roughly speaking, our idea is that $p$ and $q$ will correspond if 
$L(p \oplus q)$ is relatively small.  Notice that which projections 
then correspond to each other will strongly depend on $\rho$.  We will only consider that 
projections correspond (for a given $\rho$) if there is some
uniqueness 
to the correspondence.  The following immediate consequences of 
Proposition~\ref{prop4.3} and Theorem~\ref{th4.2} give appropriate 
expression for this uniqueness.  These consequences also hold when 
working over $\bR$.

\begin{theorem}
\label{th4.5}
Let $A$, $D$, etc., be as just above, with
$\mathrm{dist}_H^{\rho}(X,Y) < \e$.
\begin{itemize}
\item[{\em a)}] Let $p \in M_n(\cA)$ and $q \in M_n(\cD)$ be
projections, 
and suppose that $\e L(p \oplus q) < 1/2$.  If $q_1$ is any other 
projection in $M_n(\cD)$ such that $\e L(p \oplus q_1) < 1/2$ then
there 
is a path $t \mapsto q_t$ through projections in $M_n(\cD)$, going
from $q$ to $q_1$, such that
\[
L(p \oplus q_t) \le (1-\d)^{-1} \max\{L(p \oplus q),L(p \oplus q_1)\}
\]
for all $t$, where $\d = \e(L(p\oplus q) + L(p\oplus q_1))$.  If
instead there is a 
$p_1 \in M_n(\cA)$ such that $\e L(p_1 \oplus q) < 1/2$ then there is
a corresponding path from $p$ to $p_1$ with corresponding bound for 
$L(p_t \oplus q)$.
\item[{\em b)}] Let $p_0$ and $p_1$ be projections in $M_n(\cA)$ and
let 
$q_0$ and $q_1$ be projections in $M_n(\cD)$.  Set
\[
\d = \|p_0-p_1\| + \e(L(p_0 \oplus q_0) + L(p_1 \oplus q_1)).
\]
If $\d < 1$ then there are continuous paths $t \mapsto p_t$ and 
$t \mapsto q_t$ from $p_0$ to $p_1$ and $q_0$ to $q_1$, respectively, 
such that
\[
L(p_t \oplus q_t) \le (1-\d)^{-1} \max\{L(p_0 \oplus q_0),L(p_1
\oplus q_1)\}
\]
for all $t$.
\end{itemize}
\end{theorem}

We remark that a more symmetric way of stating part b) above is to 
define $\d$ by
\[
\d = \max\{\|p_0-p_1\|,\|q_0-q_1\|\} + 
\e(L(p_0 \oplus q_0),L(p_1 \oplus q_1)).
\]

Let us now examine the consequences of Corollary~\ref{cor4.4}.  
This is best phrased in terms of:

\begin{notation}
\label{not4.6}
For any $n$ let $\cP_n(X)$ denote the set of projections in
$M_n(\cA)$.  
For any $r \in \bR^+$ let
\[
\cP_n^r(X) = \{p \in \cP_n(X): L_A(p) < r\},
\]
and similarly for $Y$ and $Z$.
\end{notation}

Now $\cP_n^r(X)$ may have many path components.  We will see examples 
shortly.  As suggested in the introduction, it may be appropriate to
view 
these different path components as representing inequivalent vector
bundles, 
notably if $X$ is a finite set.  Let $\Pi$ be one of these path
components.  
Let $\Phi_X$ denote the evident restriction map from $\cP_n(Z)$ to $\cP_n(X)$ (for 
$Z = X  \dc  Y$).  For a given $s \in \bR^+$ with $s \ge r$ it may be
that 
$\Phi_X(\cP_n^s(Z)) \cap \Pi$ is non-empty.  This is an existence
question, 
which we deal with in the next sections.  But at this point, from 
Corollary \ref{cor4.4} we can conclude that:

\begin{theorem}
\label{th4.7}
Let notation be as above, with $\mathrm{dist}_H^{\rho}(X,Y) < \e$, and 
let $r \in \bR^+$ with $\e r < 1/2$.  Let $\Pi$ be a path component of 
$\cP_n^r(X)$.  Let $s \in \bR^+$ with $s \ge r$ and $\e s < 1/2$.  Let 
$p_0, \ p_1\in \Pi$ and suppose that there are $q_0$ and $q_1$ in 
$\cP_n^s(Y)$ with $L(p_j \oplus q_j) \le s$ for $j = 0,1$.  Assume, 
even more, that there is a path $\tilde p$ in $\Pi$ connecting $p_0$
and 
$p_1$ that lies in $\Phi_X(\cP_n^s(Z))$.  Then for any $\d$ with 
$2\e s < \d < 1$ there exist a path ${p}$ in $\cP_n(X)$ going from 
$p_0$ to $p_1$ and a path $q$ in $\cP_n(Y)$ going from $q_0$ to $q_1$ 
such that $L(p_t \oplus q_t) < (1-\d)^{-1}s$ for each $t$.  The
situation 
is symmetric between $X$ and $Y$, so the roles of $X$ and $Y$ can be 
interchanged in the above statement.
\end{theorem}  

Thus, in the situation described in the theorem, if $\Pi$ represents
some 
particular class of bundles on $X$, such as ``monopole'' bundles on a 
sphere, then the projections $q \in \cP_n^s(Y)$ paired with ones in
$\Pi$ 
by the requirement that $L(p \oplus q) < s$ will be homotopic, and in 
particular will determine isomorphic bundles on $Y$.  We emphasize
that 
the above pairing of projections depends strongly on $\rho$, and not
just 
on the Gromov-Hausdorff distance between $X$ and $Y$.  This reflects
the 
fact that Gromov--Hausdorff distance is only a metric on {\em isometry 
classes} of compact metric spaces.

Notice that the homotopies obtained above between $q_0$ and $q_1$ need
not 
lie in $\cP_n^s(Y)$.  We can only conclude that they lie in
$\cP_n^{s'}$ 
where $s' = (1-\d)^{-1}s$.  But at least we can say
that 
$s'$ approaches $s$ as $\e$, and so $\d$, goes to $0$.  

In Theorem~\ref{th6.4} and Corollary~\ref{cor6.7} we will deal with
the existence of actual lifts of homotopics between $p_0$ and $p_1$.


\section{Extending Lipschitz functions}
\label{sec5}

To obtain the {\em existence} of extensions to $Z$ of vector bundles
on $X$ in a manner controlled by the metric, we need to extend
projection-valued 
functions on $X$ to projection-valued functions on $Z$ with control of
the 
Lipschitz norm.  We approach this by first extending projection-valued
functions on $X$ just to general functions on $Z$ with 
values in $M_n^s(\bC)$, the space of self-adjoint matrices.  We treat
this 
problem in this section, and then in the next section we see how to
modify 
the extended functions so as to be projection-valued.

We must also treat here homotopy versions of this extension problem.  
For this purpose we let $T$ denote a compact space which will be the
parameter 
space for the homotopics, so that eventually $T$ will be an interval
in $\bR$.  We will consider functions $F$ on $X \x T$, 
and for any $t \in T$ we let $F_t$ denote the 
function $x \mapsto F(x,t)$.  When $F$ has values in a Banach space we
can then consider $L(F_t)$ for each $t$, as defined earlier. For a
discrete set $\G$ we let $\ell_{\bR}^{\i}(\G)$ denote the Banach space
of all bounded real-valued functions on $\G$ with the supremum norm.
The following proposition is a homotopy version of a well-known fact
which appears as proposition~$2.2$ of \cite{Rf76}.

\begin{proposition}
\label{prop5.1}
Let $(Z,\rho)$ be a compact metric space, and let $X$ be a closed
subset 
of $Z$.  Let $\G$ be a set (discrete, and possibly uncountable).  
Let $F$ be a continuous function from $X \x T$ to
$\ell_{\bR}^{\i}(\G)$.  
Suppose that there is a constant, $N$, such that $L(F_t) \le N$ for
all $t \in T$.  Then there is a continuous extension, $G$, of $F$ to
$Z \x T$ 
such that $L(G_t) \le N$ for all $t$, and $\|G\|_{\i} = \|F\|_{\i}$.
\end{proposition}

\begin{proof}
For each $\g \in \G$ and $x \in X$ define $H_{\g}^x$ on $Z \x T$ by
\[
H_{\g}^x(z,t) = F(x,t)(\g) - N\rho(z,x).
\]
Note that $H_{\g}^x(z,t) \le \|F\|_{\i}$.  Much as in the standard
proof of 
the non-homotopy version (as in theorem~$2.1$ of \cite{Rf76}), define 
$H: Z \x T \to \ell_{\bR}^{\i}(\G)$ by
\[
H(z,t)(\g) = \sup\{H_{\g}^x(z,t): x \in X\}.
\]
Clearly $H$ is well-defined and $H \le \|F\|_{\i}$ as functions.  It
is easily seen that for each $\g$ and $t$ 
we have $L(z \mapsto H(z,t)(\g)) \le N$, and that $H$ is an extension
of $F$.  It follows easily that $L(H_t) \le N$ for all $t \in T$.

Next, we must show that $H$ is continuous on $Z \x T$.  Let 
$(z_0,t_0) \in Z \x T$, and let $\e > 0$ be given.  Because $F$ is 
continuous and $Z \x T$ is compact, a little compactness argument
shows that there is a neighborhood, $\cN$, of $t_0$ such that for
every $x \in X$ 
and $t \in \cN$ we have $\|F(x,t) - F(x,t_0)\|_{\i} < \e/2$, and so 
$|F(x,t)(\g) - F(x,t_0)(\g)| < \e/2$ for each $\g \in \G$.  Let $\cB$
be the 
ball about $z_0$ of radius $\e/(2N)$ in $Z$.  Then for any 
$(z,t) \in \cB \x \cN$, any $x \in X$, and any $\g \in \G$ we have
\begin{eqnarray*}
|H_{\g}^x(z,t) &-& H_{\g}^x(z_0,t_0)| \\ 
&\leq& |F(x,t)(\g) - F(x,t_0)(\g)| + 
N|\rho(z,x) - \rho(z_0,x)| < \e.
\end{eqnarray*}
A simple argument very similar to the proof of proposition~$1.5.5$ of 
\cite{Wvr2} then shows that for any $(z,t) \in \cB \x \cN$ 
and $\g \in \G$ we have
\[
|H(z,t)(\g) - H(z_0,t_0)(\g)| < \e.
\]
It follows that $\|H(z,t) - H(z_0,t_0)\| \leq \e$.  Thus $H$ is
continuous.

Finally, view $-\|F\|_{\i}$ as a constant function on $Z \x T$, and
set 
$G = H \vee (-\|F\|_{\i})$, where $\vee$ means ``maximum''.  
Then $G$ has the desired properties.
\end{proof}

Now let $V$ be a finite-dimensional real Banach space (such as 
$M_n^s(\bC)$).  By definition, the projection constant, $\cP\cC(V)$,
of $V$ 
is the smallest constant $c$ such that whenever $V$ is isometrically
embedded 
into a Banach space $W$ there is a projection $P$ from $W$ onto $V$ 
such that $\|P\| \leq c$.  (Such a smallest constant exists --- see, 
e.g., proposition~$1.4$ of \cite{Rf76}.)

\begin{proposition}
\label{prop5.2}
Let $(Z,\rho)$, $X$ and $T$ be as in Proposition~$5.1$, and let $V$ be
a finite-dimensional real Banach space.  Let $F$ be a continuous
function 
from $X \x T$ into $V$ for which there is a constant, $N$, such that 
$L(F_t) \le N$ for all $t \in T$.  Then there is a continuous
extension, 
$G$, of $F$ to $Z \x T$ such that $L(G_t) \le N(\cP\cC(V))$ for all 
$t \in T$ and $\|G\|_{\i} \le \|F\|_{\i}(\cP\cC(V))$.
\end{proposition}

\begin{proof}
We can isometrically embed $V$ into $\ell_{\bR}^{\i}(\G)$ for some
discrete 
set $\G$. (For example, take $\G$ to be the unit ball of the dual
space $V'$ 
with the discrete topology.)  We can then view $F$ as a function from 
$X \x T$ into $\ell_{\bR}^{\i}(\G)$, and apply
Proposition~\ref{prop5.1} 
to find a continuous extension ${\tilde G}$ of $Z \x T$ into 
$\ell_{\bR}^{\i}(\G)$ such that $L({\tilde G}_t) \le N$ for all $t$
and 
$\|{\tilde G}\|_{\i} = \|F\|_{\i}$.  Let $P$ be a projection from 
$\ell_{\bR}^{\i}(\G)$ onto $V$ such that $\|P\| \leq \cP\cC(V)$.  
Then the function $G = P \circ {\tilde G}$ has the desired properties.
\end{proof}

For a given Banach space $V$ it is usually not easy to determine the 
precise value of $\cP\cC(V)$.  However, our projection-valued
functions 
can be viewed to have values in the Banach space $M_n^s(\bC)$ of
self-adjoint 
matrices, and in theorem~$7.2$ of \cite{Rf76} we find that
\setcounter{equation}{2}
\begin{equation}
\label{eq5.3}
\cP\cC(M_n^s(\bC)) = 2n\left( \frac {n}{n+1}\right)^{n-1} - 1.
\end{equation}
It is also noticed there that when the right-hand side is written as 
$n\o(n)$, then $\o(n)$ converges to $2e^{-1}$ as $n \to \i$.  In 
theorem~$1.5$ of \cite{Rf76} we also find that $\cP\cC(V) = \cL\cE(V)$ 
where $\cL\cE(V)$ is the Lipschitz extension constant of $V$, that is, 
the smallest constant $c$ such that every $V$-valued function $f$ on a 
subset $X$ of a metric space $Z$ can be extended to a function $g$ on
$Z$ 
such that $L(g) \le cL(f)$.  Thus in Proposition~\ref{prop5.2} above
the 
constant $\cP\cC(V)$ is the smallest constant $c$ for which we can
always find a $G$ for which $L(G) \le cL(F)$. (After \cite{Rf76} was
published I learned that the statement of theorem 1.5 of \cite{Rf76}
essentially already appears as proposition 5.1 of \cite{PrY}.) 

However, one can ask whether the constant in the inequality
$\|G\|_{\i} \le \|F\|_{\i}(\cP\cC(V))$ of Proposition~\ref{prop5.2}
can be improved.  In some cases, such as when $V$ is a Hilbert space,
$\cP\cC(V)$ can be replaced by $1$ in this inequality.  But I have not
been able to determine whether this constant can be improved for the
case of $V = M_n^s(\bC)$.  However, in proposition~$8.3$ of
\cite{Rf76} it is shown by means of radial retractions that for any
$V$ this constant can be replaced by $1$ if the Lipschitz inequality
is weakened to $L(q) \le L(f)(2\cP\cC(V))$.  This applies equally well
to the homotopy situation we consider here.  We will use this fact in
the rest of our paper, since to have $\|\cdot\|_{\i}$ preserved under
extensions will simplify our bookkeeping. If eventually situations are
found where this has undesirable effects, the extra bookkeeping can be
done easily.  For the rest of this paper we will use:

\setcounter{theorem}{3}
\begin{notation}
\label{not5.4}
For each positive integer $n$ we let $\l_n$ denote the smallest
constant $c$ such that, with notation as above, for any $Z$, $X$, $T$,
and for any continuous $F: X \x T \to M_n^s(\bC)$ for which there is a
constant $N$ such that $L(F_t) \le N$ for all $t$, there is a
continuous extension, $G$, of $F$ to $Z \x T$ such that
\[
L(G_t) \le \l_nN \mbox{ for all $t$ 
\quad and \quad $\|G\|_{\i} = \|F\|_{\i}$}.
\]
\end{notation}

Thus if we let $\l_n^*$ be given by the formula in \eqref{eq5.3}, then
we have
\[
\l_n^* \le \l_n \le 2\l_n^*.
\]

Suppose we are instead working over $\bR$.  Given a function $F: X \x
T \to M_n^s(\bR)$, we can view it as taking values in $M_n^s(\bC)$,
and then find an extension, ${\tilde G}$, of it to $Z \x T$ which
satisfies the estimates of Notation~\ref{not5.4}.  But the map which
for each matrix in $M_n^s(\bC)$ replaces all of its entries by their
real part is norm-decreasing and $\bR$-linear.  Let $G$ be the
composition of this map with ${\tilde G}$.  We thus obtain:

\begin{proposition}
\label{prop5.5}
Let $F$ be a continuous function from $X \x T$ to $M_n^s(\bR)$ for
which there is a constant $N$ such that $L(F_t) \le N$ for each $t$.
Then there is a continuous extension, $G: Z \x T \to M_n^s(\bR)$, of
$F$ such that $L(G_t) \le \l_nN$ for each $t$, and $\|G\|_{\i} =
\|F\|_{\i}$.
\end{proposition}

One can probably replace $\l_n$ by a smaller constant by using the
techniques of section~$7$ of \cite{Rf76} to compute
$\cP\cC(M_n^s(\bR))$.

There is another aspect of extensions which might seem relevant,
namely that one can define certain classes of compact metric spaces
for which the constant for extending Lipschitz functions into any
Banach space is smaller than for arbitrary compact metric spaces.  For
example, in \cite{JLS} it is shown that there is a universal constant
$C$ such that if $Z$ is a compact subset of some $d$-dimensional
Banach space, with metric from the Banach space, then for every closed
subset $X$ of $Z$ and for every Banach space $V$, every function $f$
from $X$ to $V$ can be extended to a function $g$ from $Z$ to $V$ such
that $L(q) \le CdL(f)$.  Notice that this inequality is independent of
the dimension of $V$, unlike our results above.  Even more, from
part~$5$ of theorem~$5.1$ of \cite{LN} we see that if $Z$ is in  fact
a compact subset of a $d$-dimensional Hilbert space, then the above
inequality can be improved to $L(q) \le C(d^{1/2})L(f)$. See also
\cite{BrB}.

But these results do not seem to be useful to us here for the
following reason.  Ultimately we want to consider two compact metric
spaces $(X,\rho_X)$ and $(Y,\rho_Y)$, and then, in order to use their
Gromov--Hausdorff distance we must consider all isometric embeddings
of them into compact metric spaces $(Z,\rho)$.  If we were to restrict
$X$ and $Y$ to be subsets of a $d$-dimensional Banach space, I see no
reason why we could, for example, restrict $(Z,\rho_Z)$ to also be
(isometrically) a subset of a $d$-dimensional Banach space.  But this
type of question would be interesting to explore.


\section{Extending vector bundles}
\label{sec6}

We extend vector bundles in a controlled way by extending their
projections.  As before, let $\pi$ denote the restriction map from $B
= C(Z)$ onto $A = C(X)$, and let $p$ be a projection in $M_n(\cA)$.
In terms of Notation~\ref{not5.4} we can find $b \in M_n^s(B)$ such
that $\pi(b) = p$, $L(b) \le \l_nL(p)$ and $\|b\| = \|p\| = 1$.  (The
fact that $\|b\| = 1$ is one of the places where our definition of
$\l_n$ simplifies the bookkeeping.)  Suppose now that $X$ is
$\e$-dense in $Z$.  Then from Key Lemma~\ref{keylem4.1} we see that
\begin{eqnarray*}
\|b^2-b\| &\le& \|\pi(b^2-b)\| + \e L(b^2-b) \\ 
&\le& \|p^2-p\| + \e(2\|b\|L(b) + L(b)) = \e 3L(b).
\end{eqnarray*}
Note the crucial use made here of the Leibniz property of $L$.  Set
$\d = \e 3L(b)$, so that $\|b^2-b\| \le \d$.  Much as in
Lemma~\ref{lem3.2} we have:

\begin{lemma}
\label{lem6.1}
Let $b$ be an element in a unital $C^*$-algebra such that $b^* = b$
and $\|b^2-b\| \le \d$.  Then
\[
\s(b) \subseteq [-2\d,2\d] \cup [1-2\d,1+2\d].
\]
\end{lemma}

\begin{proof}
Let $\l \in \s(b)$.  The polynomial $(x^2-x) - (\l^2-\l)$ has value
$0$ at $x=\l$, and so factor as $(x-\l)p(x)$ for some polynomial $p$.
On substituting $b$ for $x$ in this factorization, one sees easily
that $\l^2 - \l \in \s(b^2-b)$.  Thus $|\l^2-\l| \le \d$, that is,
$|\l| |\l-1| \le \d$.  If $|\l| \ge 1/2$, 
then $(1/2)|\l-1| \le |\l| |\l-1| \le \d$, so that $|\l-1| \le 2\d$
and $\l \in [1-2\d,1+2\d]$.  If $|\l| \le 1/2$ then $|\l-1| \ge 1/2$,
so that $|\l|(1/2) \le |\l| |\l-1| \le \d$.  Thus $|\l| \le 2\d$ so
that $\l \in [-2\d,2\d]$.
\end{proof}

Suppose now that $\d < 1/4$, so that the two intervals $[-2\d,2\d]$ and
$[1-2\d,1+2\d]$ are disjoint.  Then we are in almost the same
situation as in the proof of Proposition~\ref{prop3.3}, but with
slightly different conditions on the spectrum.  We now present the
result that we seek, but in the greater generality involving
homotopies, much as in Notation~\ref{not5.4}.  In particular, let
$\l_n$ be as in Notation~\ref{not5.4}, let $A = C(X)$, etc., and let
$T$ be a compact space which serves as a parameter space.

\begin{theorem}
\label{th6.2}
Let $p: T \to M_n(A)$ be a continuous function such that $p_t$ is a
projection for each $t$.  Assume that there is a constant, $N$, such
that $L(p_t) \le N$ for all $t$.  Suppose that $X$ is $\e$-dense in
$Z$.  If $\e \l_n N < 1/12$, then there exists a continuous function $q: T
\to M_n(B)$ such that $q_t$ is a projection, $\pi(q_t) = p_t$, and 
\[
L(q_t) <  \l_n N(1 - 12\e \l_nN)^{-1},
\]
for each $t \in T$.  If, instead, $A = C_{\bR}(X)$, etc., one has the
same conclusions.
\end{theorem}

\begin{proof}
We can view $p$ as a continuous function on $X \times T$. Then
according to the definition of $\l_n$ in Notation~\ref{not5.4} we can
find a continuous function $G: T \to M_n(B)$ such that for 
each $t \in T$ we have 
$G_t^* = G_t$, $\pi(G_t) = p_t$, $\|G_t\| = \|p_t\| = 1$ 
and $L(G_t) \le \l_n L(p_t) \le \l_n N$.  From the discussion given
before Lemma~\ref{lem6.1} we see that for any $t \in T$ we have
$\|G_t^2 - G_t\| \le \e 3L(G_t) \le \e 3\l_n N$.  
Let $\d = \e 3\l_nN$.  Then according to Lemma~\ref{lem6.1} and the
comment immediately after it, if $\d < 1/4$ then $\s(G_t)$ is
contained in the union of the disjoint intervals $[-2\d,2\d]$ and
$[1-2\d,1+2\d]$.  So we now assume that $\d < 1/4$.  Much as in the
proof of Proposition~\ref{prop3.3}, let $\chi$ be defined on $\bC$ by
setting $\chi(z) = 0$ if $\mathrm{Re}(z) \le \d + \frac {1}{4}$, and
$\chi(z) = 1$ otherwise.   For any nice curve $\g$ around
$[1-2\d,1+2\d]$ that lies in the domain where $\chi$ is
holomorphic we now set, for each $t \in T$,
\[
q_t = \frac {1}{2\pi i} \int_{\g} \chi(z)(z - G_t)^{-1}dz.
\]
Then as in the proof of Proposition~\ref{prop3.1} we see that $q_t$ is
a projection in $M_n(\cB)$.  From the facts that $\pi(G_t) = p_t$ and
that $p_t$ is a projection it is easily verified 
that $\pi(q_t) = p_t$.

We need to estimate $L(q_t)$. Instead of using the curve that we used
in earlier versions of this paper we can now directly apply proposition 3.1 
of \cite{Li} much as we did in the proof of
Proposition~\ref{prop3.3}.  This gives
\[
L(q_t) \leq L(G_t)(1-4\d)^{-1} \leq \l_n N(1-12\e \l_nN)^{-1}.
\]

Finally, we must show that $q$ is continuous in $t$.  
Given $s,t \in T$ we have, by a familiar maneuver (basically the
``resolvant equation'' \cite{Kto}),
\begin{eqnarray*}
q_t - q_s &= &\frac {1}{2\pi i} \int_{\g} \chi(z)((z-G_t)^{-1} -
(z-G_s)^{-1})dz \\
&= &\frac {1}{2\pi i} \int_{\g} (z-G_t)^{-1}((z-G_s) -
(z-G_t))(z-G_s)^{-1}dz \\
&= &\frac {1}{2\pi i} \int_{\g} (z-G_t)^{-1}(G_t-G_s)(z-G_s)^{-1} dz.
\end{eqnarray*}
Thus
\[
\|q_t-q_s\| \le \|G_t-G_s\| \frac {1}{2\pi} \int_{\g} \|(z-G_t)^{-1}\|
\|(z-G_s)^{-1}\| d|z| \le K\|G_t - G_s\|
\]
for a suitable constant $K$ obtained by the kind of estimates used
above to bound $L(q_t)$.  Since $G$ is continuous, it follows that $q$
is also. 

On looking at Proposition~\ref{prop2.4} and its proof it is easy to
see how to adapt the above proof to the case in which $A =
C_{\bR}(X)$, etc. 
\end{proof}

\begin{corollary}
\label{cor6.3}
Let $p$ be a projection in $M_n(A)$.  Assume that $X$ is $\e$-dense in
$Z$.  If $\e \l_nL(p) < 1/12$, 
then there exists a projection $q \in M_n(B)$ such that $\pi(q) = p$ and
\[
L(q) < \l_n L(p)(1-12\e \l_nL(p))^{-1}.
\]
\end{corollary}

We remark that the role of the Cauchy integrals used in this section can
be viewed as follows: The set of elements $b \in  M^s_n(A)$ such that
$\|b^2 -b\| < 1/4$ is a neighborhood of the set of projections in  $M^s_n(A)$,
and the Cauchy integrals give a retraction from this neighborhood onto the set
of projections in such a way that one keeps control of the Lipschitz constants.

It would be interesting to know whether different techniques, perhaps
not involving at all the Lipschitz extension properties of functions
into $M_n^s(\bC)$ as used in Section~\ref{sec5}, but rather maybe
working directly just with functions into the space of projections
(or into a Grassman manifold if $X$ is connected), for example
by using in part the methods of \cite{Pds},  
could yield extensions with smaller increase in the Lipschitz
constants than is obtained in the above theorem and corollary.

Notice that Corollary~\ref{cor6.3} gives a criterion for
extending a vector bundle from $X$ to $Z$ which is quite independent
of how complicated the topologies of $X$ and $Z$ are.  All that is
required is that there be a metric on $Z$ and a projection $p \in
M_n(C(X))$ representing the bundle such that $X$ is $\e$-dense in $Z$
and $\e\l_n L(p) < 1/12$.  There is no requirement that any
obstructions from algebraic topology vanish, or that spaces have
finite dimension or be locally geometrically $n$-connected
as seems to be needed in \cite{Pds}.

We now combine Theorem~\ref{th6.2} with the uniqueness given by
Proposition~\ref{prop4.3} to treat the case of a path from $p_0$ to
$p_1$ for which we already have lifts of $p_0$ and $p_1$.

\begin{theorem}
\label{th6.4}
Let $p: [0,1] \to M_n(A)$ be a path of projections, and let $q_0$ and
$q_1$ be projections in $M_n(B)$ such that $\pi(q_0) = p_0$ and
$\pi(q_1) = p_1$.  Suppose that there is a constant, $N$, such that
$L(p_t) \le N$ for all $t$.  Set $N' = \max\{L(q_0),L(q_1)\}$.  Assume
that $X$ is $\e$-dense in $Z$.  
If $\e \l_n N < 1/14$ and $\e N' < 1/2$ then there exists a path, $q$, of
projections from $q_0$ to $q_1$ such that $\pi(q_t)$ is in the range
of the path $p$ for each $t$ (though $\pi(q)$ may have a different
parametrization) and 
\[
L(q_t) \le ((1/2) - \e N')^{-1} \max\{7\l_n N,N'\}.
\]
The same conclusion holds if $A = C_{\bR}(X)$, etc.
\end{theorem}

\begin{proof}
From Theorem~\ref{th6.2} we see that there is a path ${\tilde q}$ of
projections in $M_n(B)$ such that $\pi({\tilde q}_t) = p_t$ and
\[
L({\tilde q}_t) \le \l_n N(1-12\e \l_n N)^{-1} \le \l_n N(1-(6/7))^{-1} =
7\l_n N
\]
for all $t \in [0,1]$.  In particular, for all $t$
\[
\e L({\tilde q}_t) \le 7\e \l_n N  < 7/14 = 1/2.
\]
Since $\pi({\tilde q}_0) = p_0 = \pi(q_0)$ and $\e L(q_0) \le \e N' <
1/2$, we can apply Proposition~\ref{prop4.3} to obtain a path of
projections, $t \mapsto q_t^0$, joining $q_0$ to ${\tilde q}_0$ and
such that for each $t$ we have $\pi(q_t^0) = p_0$ and
\begin{eqnarray*}
L(q_t^0) &\le &(1-\e(L(q_0) + L({\tilde q}_0)))^{-1}
\max\{L(q_0),L({\tilde q}_0)\} \\
&\le &(1-\e L(q_0) - 1/2)^{-1} \max\{L(q_0), 7\l_n N\} \\
&\le &((1/2)-\e N')^{-1} \max\{N',7\l_nN\}.
\end{eqnarray*}
In the same way there is a path of projections, $t \mapsto q_t^1$,
connecting ${\tilde q}_1$ to $q_1$ with corresponding bound on
$L(q_t^1)$.  We concatenate the three paths $q^0$, ${\tilde q}$ and
$q^1$ to obtain a path, $q$, of projections connecting $q_0$ to $q_1$
such that each $\pi(q_t)$ is in the range of the path $p$.

Since $((1/2)-\e N')^{-1} > 1$ and $L({\tilde q}_t) < 7\l_n N$, we
see that the bound given above for $L(q_t^0)$ and $L(q_t^1)$ is also a
bound for $L({\tilde q}_t)$, and thus for $L(q_t)$ for all $t$.
\end{proof} 

Let us now see what consequences the above existence results have for
metric spaces that are close together.  As we did near the end of
Section~\ref{sec4}, we let $Z = X  \dc  Y$, with $\rho$ on $Z$
restricting to the given metrics on $X$ and $Y$.  We also let 
$D = C(Y)$ as before, so that $B = A \oplus D$.  We use the notation
$\cP^r(X)$, etc., introduced in Notation \ref{not4.6}.  By using
Corollary~\ref{cor6.3} to satisfy the hypothesis concerning
$\Phi_X(\cP_n^s(Z))$ in Theorem~\ref{th4.7} we obtain:

\begin{theorem}
\label{th6.5}
Let $r \in \bR^+$ be given.  Let $\e$ be small enough that 
$\e \l_n r < 1/14$.  Set $s = \l_n r(1 - 12\e \l_n r)^{-1}$, 
so that $\e s < 1/2$.  Finally,
assume that $\mathrm{dist}_H^{\rho}(X,Y) < \e$.  Let $p_0$ and $p_1$
be projections in $\cP_n^r(X)$ which lie in the same path component of
$\cP_n^r(X)$.  By Corollary~$\ref{cor6.3}$ there exist $q_0$ and $q_1
\in \cP^s(Y)$ such that $L(p_j \oplus q_j) < s$ for $j = 0,1$.  For
any such $q_0$ and $q_1$ and for any $\d$ with $2\e s < \d < 1$ there
exist a path $p$ in $\cP_n(X)$ going from $p_0$ to $p_1$ and a path
$q$ in $\cP_n(Y)$ going from $q_0$ to $q_1$ such that
\[
L(p_t \oplus q_t) < (1-\d)^{-1}s
\]
for all $t$.  In particular, the vector bundles determined by $q_0$
and $q_1$ are isomorphic.
\end{theorem}

The important comments made in the two paragraphs following
Theorem~\ref{th4.7} apply equally well to Theorem~\ref{th6.5}.

Suppose that we have a specific homotopy in $\cP_n(X)$ and we want a
homotopy in $\cP_n(Y)$ which corresponds to it, but we do not have
specific endpoints in $\cP_n(Y)$ that we require be joined by the
homotopy.  Then we can apply Theorem~\ref{th6.2} to obtain:

\begin{theorem}
\label{th6.6}
Let $r \in \bR^+$ be given.  Let $\e$ be small enough that 
$\e \l_n r < 1/12$.  Assume that $\mathrm{dist}_H^{\rho}(X,Y) < \e$.
Set $s = \l_n r(1-12\e \l_n r)^{-1}$.  Then for any path $p$ in
$\cP_n^r(X)$ there exists a path $q$ in $\cP_n^s(Y)$ such that 
$L(p_t \oplus q_t) < s$ for every $t$.
\end{theorem}

Suppose finally that we have a specific homotopy $p$ in $\cP_n(X)$ and
specific $q_0$ and $q_1$ in $\cP_n(Y)$ which correspond to $p_0$ and
$p_1$ for $\rho$, and we want a corresponding homotopy in $\cP_n(Y)$
joining $q_0$ and $q_1$.  We can apply Theorem~\ref{th6.4} to
immediately obtain:

\begin{corollary}
\label{cor6.7}
Let $r \in \bR^+$ be given and let $\e$ be small enough that 
$\e \l_n r < 1/14$.  Let $p$ be a path in $\cP_n^r(X)$.  Let $q_0$ and
$q_1$ be projections in $M_n(\cD)$, and set 
$N = \max\{L(p_0 \oplus q_0),L(p_1 \oplus q_1)\}$.  Assume further
that $\e$ is small enough that $\e N < 1/2$.  Finally, assume that
$\mathrm{dist}^\rho_H(X,Y) < \e$.  Then there exists a path $q$ of
projections in $M_n(\cD)$ going from $q_0$ to $q_1$, and a
reparametrization ${\tilde p}$ of the path $p$, still with domain
$[0,1]$, such that
\[
L({\tilde p}_t \oplus q_t) \le ((1/2)-\e N)^{-1} \max\{7\l_nr,N\}
\]
for all $t$.
\end{corollary}

Of course, vector bundles can be represented by projections of
different sizes.  In particular, if $p \in \cP_n(X)$, then for $m > n$
the projection 
${\tilde p} = \begin{pmatrix} p & 0 \\ 0 & 0 \end{pmatrix}$ in
$\cP_m(X)$, for the $0$'s of appropriate size, will represent the same
bundle as $p$, and we will have $L({\tilde p}) = L(p)$.  But because
$\l_n$ grows with $n$, I have so far not seen anything really useful
to say about how projections of different sizes should be related
within our context of Gromov--Hausdorff distance.


\section{Projective modules and frames}
\label{sec7}

We now make some preparations for our discussion of specific examples.
Naturally-arising vector bundles are not often presented by means of
projections, and there is usually no canonical choice of a projection
for them.  We recall in this section some elementary tools for
obtaining projections corresponding to vector bundles.

As before, we set $A = C(X)$ for $X$ a compact space.  (With evident
modifications, everything in this section works just as well for $A =
C_{\bR}(X)$.)  Most of the discussion in this section applies without
change to a general unital $C^*$-algebra $A$, and so we will in some
places write it in that generality, but the reader can take $A$ to be
$C(X)$ with no disadvantage for reading the next sections.  Let $\Xi$
be an $A$-module.  We use right-module notation, both because it eases
the bookkeeping somewhat, and also in view of the generalizations that
we will consider elsewhere in which $A$ is non-commutative, for which
most writers use right modules.  By an $A$-valued inner product on
$\Xi$ (for example, a Riemannian or Hermitian metric on $\Xi$
according to whether we work over $\bR$ on $\bC$) we mean \cite{RLL,
WO} a sesquilinear form $\< \cdot,\cdot \>_A$ on $\Xi$ with values in
$A$ such that for $\xi,\eta \in \Xi$ and $a \in A$ we have
\begin{itemize}
\item[1)] $\<\xi,\eta a\>_A = \<\xi,\eta\>_A a$,
\item[2)] $(\<\xi,\eta\>_A)^* = \<\eta,\xi\>_A$ (with $* =$ complex
conjugation),
\item[3)] $\<\xi,\xi\>_A \ge 0$, with $= 0$ only if $\xi = 0$.
\end{itemize}
For naturally-arising vector bundles there is often a natural choice
of $C(X)$-valued inner product, even when there is no natural choice
of projection.  We will see this in the examples in the later
sections.

On $A^n$ as a right $A$-module we have the standard inner product
defined by
\[
\<(a_j),(b_k)\>_A = \sum a_j^*b_j.
\]

If $\Xi = pA^n$ for some projection $p \in M_n(A)$, then the
restriction to $\Xi$ of the inner product on $A^n$ will be an inner
product on $\Xi$.  Thus every (finitely generated) projective
$A$-module (that is, a summand of $A^n$ for some $n$) has an inner product.  If we set $\eta_j = pe_j$ for each
$j$, where $\{e_j\}$ is the ``standard basis'' for $A^n$, then
$\{\eta_j\}$ is a ``standard module frame'' for $\Xi$.  We recall
\cite{FrL1, FrL2, Rf69} the general definition, valid for modules over
any unital $C^*$-algebra (over $\bC$ or $\bR$). 

\begin{definition}
\label{def7.1}
Let $A$ be a unital $C^*$-algebra and let $\Xi$ be a right $A$-module.
Let $\Xi$ be equipped with an $A$-valued inner-product,
$\<\cdot,\cdot\>_A$.  By a (finite) {\em standard module frame} for
$\Xi$ (with respect to the inner-product) we mean a finite family
$\{\eta_j\}$ of elements of $\Xi$ such that for any $\xi \in \Xi$ the
reconstruction formula
\[
\xi = \sum \eta_j\<\eta_j,\xi\>_A
\]
is valid.
\end{definition}

The relationships that we need between standard module frames and the
projections corresponding to projective modules are given (see, for
example, scattered places in \cite{Rf33, FrL1, FrL2}) by:

\begin{proposition}
\label{prop7.2}
Let $\Xi$ be a right module over a unital $C^*$-algebra, $A$, and
suppose that $\Xi$ is equipped with an $A$-valued inner-product.  If
$\Xi$ has a standard module frame, $\{\eta_j\}_{j=1}^n$, then $\Xi$ is
a projective $A$-module.  In fact, $\Xi \cong pA^n$ isometrically,
where $p$ is the projection in $M_n(A)$ defined by 
$p_{jk} = \<\eta_j,\eta_k\>_A$.  Furthermore, $\Xi$ is self-dual for
its inner product, in the sense that for any $\var \in
\mathrm{Hom}_A(\Xi,A_A)$ there is a (unique) \ $\z_{\var} \in \Xi$
such that 
$\var(\xi) = \<\z_\var,\xi\>_A$ for all $\xi \in \Xi$.
\end{proposition}

\begin{proof}
Let $\{\eta_j\}^n_{j=1}$ be a standard module frame for $\Xi$.  Define
$\Phi: \Xi \to A^n$ by
\[
(\Phi\xi)_j = \<\eta_j,\xi\>_A.
\]
Clearly $\Phi$ is an $A$-module homomorphism.  From the reconstruction
formula in the definition of a standard module frame it is clear that
$\Phi$ is injective.  Clearly $p^* = p$.  Furthermore
\begin{eqnarray*}
(p^2)_{ik} &= &\sum_j p_{ij}p_{jk} = \sum_j \<\eta_i,\eta_j\>_A \<
\eta_j,\eta_k\>_A \\
&= &\left\<\eta_i, \ \sum \eta_j\<\eta_j,  \eta_k\>_A\right\>_A =
\<\eta_i,\eta_k\>_A = p_{ik}.
\end{eqnarray*}
Thus $p^2 = p$, and so $p$ is a projection.  Now for every $\xi \in
\Xi$ we have
\begin{eqnarray*}
(p(\Phi\xi))_j &= &\sum \<\eta_j,\eta_k\>_A(\Phi\xi)_k = \sum
\<\eta_j,\eta_k\>_A\<\eta_k,\xi\>_A \\
&= &\left\<\eta_j, \sum \eta_k\<\eta_k,\xi\>_A\right\>_A =
\<\eta_j,\xi\>_A = (\Phi\xi)_j.
\end{eqnarray*}
Thus the range of $p$ contains the range of $\Phi$.  On the other hand
if $v$ is an element of $A^n$ in the range of $p$, so that $v = pv$,
then, since $v_j \in A$,
\[
v_j = \sum\<\eta_j,\eta_k\>_Av_k = \left\<\eta_j,\sum
\eta_kv_k\right\>_A.
\]
Thus if we set $\xi = \sum \eta_kv_k$, then $v = \Phi\xi$.  Hence $p$
is exactly the projection onto the range of $\Phi$.  It follows that
$\Xi$ is a projective module.

We now show that $\Phi$ is isometric.  For $\xi,\z \in \Xi$ we have
\begin{eqnarray*}
\<\Phi\xi,\Phi\z\>_A &= &\sum \<\eta_j,\xi\>_A^*\<\eta_j,\z\>_A \\
&= &\left\<\xi,\sum \eta_j\<\eta_j,\z\>_A\right\>_A = \<\xi,\z\>_A.
\end{eqnarray*}

Finally, we show that $\Xi$ is self-dual for its inner product.  Let
$\var \in \mathrm{Hom}_A(\Xi,A_A)$, where $A_A$ means that $A$ is
viewed as a right module over itself.  Then for any $\xi \in \Xi$ we
have
\[
\var(\xi) = \var\left( \sum \eta_j\<\eta_j,\xi\>_A\right) = \sum
\var(\eta_j)\<\eta_j,\xi\>_A = \left\<\sum
\eta_j(\var(\eta_j))^*,\xi\right\>_A.
\]
Thus $\z_{\var} = \sum \eta_j(\var(\eta_j))^*$ is the desired element
of $\Xi$ representing $\var$.
\end{proof}

Suppose that we have a projective module $\Xi$ that is already
equipped with an $A$-valued inner product $\<\cdot,\cdot\>_A$ for
which it is self-dual, and suppose that $\Phi$ is an isomorphism from
$\Xi$ to $pA^n$ for some projection $p$.  Let $\<\cdot,\cdot\>'_A$
denote the pull-back to $\Xi$ of the restriction to $pA^n$ of the
standard inner product on $A^n$.  From the self-duality of the inner
products it is easily seen that there is an $S \in
\mathrm{End}_A(\Xi)$ which is invertible and positive (for either
inner product) such that
\[
\<\xi,\eta\>_A = \<S\xi,S\eta\>'_A
\]
for all $\xi,\eta \in \Xi$.  Define $\Phi': \Xi \to pA^n$ by 
$\Phi'(\xi) = \Phi(S\xi)$.  Then for $\xi,\eta \in \Xi$ we have
\[
\<\Phi'\xi,\Phi'\eta\>_A = \<\Phi(S\xi),\Phi(S\eta)\>_A = 
\<S\xi,S\eta\>'_A = \<\xi,\eta\>_A.
\]
Thus for a given self-dual inner product on $\Xi$, and for any
projection 
$p$ representing $\Xi$ we can assume that our isomorphism 
$\Phi: \Xi \to pA^n$ preserves the inner products (i.e., is
``isometric'').  
Then on setting $\eta_j = \Phi^{-1}(pe_j)$ for each $j$ we obtain a 
standard module frame for $\Xi$ such that
\[
p_{jk} = \<\eta_j,\eta_k\>_A
\]
for all $j,k$.  We thus obtain:

\begin{proposition}
\label{prop7.3}
Let $\Xi$ be a projective $A$-module equipped with a fixed $A$-valued 
inner product for which it is self-dual.  Every projection $p$ such 
that $\Xi \cong pA^n$ for some $n$ is of the form
\[
p_{jk} = \<\eta_j,\eta_k\>_A
\]
for some standard module frame $\{\eta_j\}$ for $\Xi$.
\end{proposition}

Thus, in the presence of a metric $\rho$ on $X$, to calculate $L(p)$
for various projections $p$ representing $\Xi$ it suffices to consider
standard module frames and their corresponding projections.


\section{The M\"obius strip}
\label{sec8}

In this section and the next we show that the simplest non-trivial
vector bundle, the M\"obius-strip bundle, already provides interesting
examples that illustrate our general theory.  This requires working
over $\bR$.

Let $\bT$ denote the circle, viewed either as $\bR/\bZ$, or as 
$I = [0,1]$ with endpoints identified.  Let $A = C_{\bR}(\bT)$, which
we will usually view as consisting of functions on $\bR$ periodic of
period~$1$.  As before, we equip the free $A$-module $A^n$ with its
``standard'' inner-product, defined by
\[
\<v,w\>_A(r) = \sum v_j(r)w_j(r)
\]
for $v,w \in A^n$.  We take as our metric $\rho$ the metric coming
from the absolute-value on $\bR$.  Thus
\[
\rho(r,s) = \min\{|r-s-n|: n \in \bZ\}.
\]
We let $L$ denote the corresponding Lipschitz seminorm on $A$.

\begin{notation}
\label{not8.1}
The M\"obius-strip $A$-module $\Xi$ consists of the $\bR$-valued
continuous functions $\xi$ on $\bR$ which satisfy the condition that
for any $r \in \bR$
\[
\xi(r-1) = -\xi(r).
\]
The action of $A$ on $\Xi$ is by pointwise multiplication of
functions.  We define an $A$-valued inner-product (Riemannian metric)
on $\Xi$ by
\[
\<\xi,\eta\>_A(r) = \xi(r)\eta(r).
\]
\end{notation}

If $\Xi$ were a free $A$-module then it would contain an element $\xi$
such that $\<\xi,\xi\>_A$ is nowhere $0$, which is easily seen not to
happen.  So we seek standard module frames for $\Xi$.  Suppose that
$\{\eta_j\}_{j=1}^n$ is a standard module frame for $\Xi$, and define
a function $u: \bR \to \bR^n$ by $u(r) = (\eta_j(r))_{j=1}^n$.  It is
easily seen that $\mathrm{End}_A(\Xi)$ can be identified with $A$
itself (essentially because $\Xi$ comes from a line bundle).  The
reconstruction formula for $\{\eta_j\}$ then implies that $\|u(r)\|=1$
for all $r$, where the norm here is the Euclidean norm on $\bR^n$.
Because $\eta_j \in \Xi$ for each $j$, we also have $u(r-1) = -u(r)$
for each $r$.  It is easily seen that conversely, if $u$ is a
continuous function from $\bR$ to $\bR^n$ such that $\|u(r)\| = 1$ and
$u(r-1) = -u(r)$ for each $r$, then the component functions of $u$
form a standard module frame for $\Xi$.  For a standard module frame
$\{\eta_j\}$ and its $u$, the corresponding projection $p$ has as
entries $p_{jk}(r) = \eta_j(r)\eta_k(r)$ at $r$.  From this we easily
see that $p(r)$ is just the rank-$1$ projection onto $u(r)$, which we
like to denote by $\<u(r),u(r)\>_0$.  Briefly, $p = \<u,u\>_0$.
(Notice that $p(r-1) = p(r)$.)  

For now and later we need the undoubtedly well-known:

\begin{proposition}
\label{prop8.2}
Let $\cH$ be a Hilbert space over $\bR$ or $\bC$, and let 
$v,w \in \cH$ with $\|v\| = 1 = \|w\|$, and with corresponding
rank-$1$ projections $\<v,v\>_0$ and $\<w,w\>_0$.  Then
\[
\|\<v,v\>_0 - \<w,w\>_0\| = (1 - |\<v,w\>_{\cH}|^2)^{1/2} \le \|v-w\|.
\]
If $\cH$ is over $\bR$, the middle term is equal to $|\sin \th|$ where
$\th$ is the angle between $v$ and $w$.
\end{proposition}

\begin{proof}
If $w = av$ with $a \in \bC$ and $|a| = 1$ then the left-hand side is
$0$.  If $v$ and $w$ are linearly independent, let $\{e_1,e_2\}$ be an
orthonormal basis for the subspace spanned by $v$ and $w$, with 
$e_1 = v$.  Let $w = ae_1 + be_2$ for scalars $a$ and $b$, so that
$|a|^2 + |b|^2 = 1$.  Let $T = \<v,v\>_0 - \<w,w\>_0$.  Then the
matrix for $T$ for the basis $\{e_1,e_2\}$ is
\[
\begin{pmatrix}
1 - a{\bar a} & -a{\bar b} \\
-{\bar a}b & -b{\bar b}
\end{pmatrix} = \begin{pmatrix}
b{\bar b} & -a{\bar b} \\
-{\bar a}b & -b{\bar b}
\end{pmatrix} .
\]
Its trace is $0$ and its determinant  is $-|b|^2$.  Thus its norm is
$|b| = (1 - |\<v,w\>_{\cH}|^2)^{1/2}$, while
\[
\|v-w\|^2 = |1 - a|^2 + |b|^2 \ge |b|^2,
\]
giving the desired inequality.
\end{proof}

Because for our examples most of our spaces will be manifolds and we
will use the letters $X,Y$ for vector fields, we will at times denote
our metric space by $M$.

\begin{corollary}
\label{cor8.3}
Let $(M,\rho)$ be a metric space and let $\cH$ be a Hilbert space
(over $\bR$ or $\bC$).  Let $u$ be a function from $M$ to $\cH$ with
$\|u(m)\| = 1$ for all $m \in M$.  Define a function $p$ by 
$p(m) = \<u(m),u(m)\>_0$ for all $m \in M$.  Then $L(p) \le L(u)$.
\end{corollary}

\begin{proof}
For $m,n \in M$ we have from Proposition~\ref{prop8.2} $\|p(m) -
p(n)\| \le \|u(m) - u(n)\|$.  Now divide by $\rho(m,n)$.
\end{proof}

Since $\bT$ is a manifold, it is helpful to use calculus.  
Because we have chosen a metric that is invariant under ``rotation''
of $\bT$, we can apply Proposition~\ref{prop2.5}.  Specifically, we
will apply that proposition to functions on $\bR$ which are periodic
of period two, and so to the components of a function 
$u: \bR \to \bR^n$ that satisfies $\|u(r)\| = 1$ and $u(r-1) = -u(r)$
for each $r \in \bR$.  We conclude that for any $\e > 0$ such a
function can be approximated by a function $h: \bR \to \bR^n$ which is
infinitely differentiable, in such a way that $\|u-h\|_{\i} < \e$ and
$L(h) \le L(u)$.  Furthermore, the smoothing argument in the proof of
Proposition~\ref{prop2.5} can easily be seen to give $h(r-1) = -h(r)$.
But we need to obtain a smooth  unit-vector-valued function in order
to obtain a projection.  From the relation to $u$ we see that
$\|h(r)\| \ge 1-\e$ for all $r$.  Define $g$ by 
$g(r) = (1-2\e)^{-1}h(r)$.  Then $\|g(r)\| > 1$ for all $r$, and 
$L(g) \le (1-2\e)^{-1}L(u)$.  Also, 
$\|h-g\| \le 2\e(1-2\e)^{-1}\|h\| \le 2\e(1+\e)(1-2\e)^{-1}$, with
corresponding estimate for $\|u-g\|$.  Finally, let $v$ be the
composition of $g$ with the radial retraction from $\bR^m$ onto its
unit ball.  The radial retraction for a Hilbert space has Lipschitz
constant $1$, and is smooth at points strictly outside the unit ball.
It follows that $v$ is smooth, that $L(v) \le (1-2\e)^{-1}L(u)$, and
that $v$ can be as close to $u$ as desired by making $\e$ small
enough.  Let $q = \<v,v\>_0$.  Then $q$ is smooth, 
$L(q) \le (1-2\e)^{-1}L(p)$, and $q$ can be as close to $p$ as
desired. (Of course, we could have applied Theorem \ref{th3.5}
here.)

The above arguments can also be used for some other examples involving
real line bundles over certain manifolds.  We state the next step in
somewhat general form in order to contrast the situation over $\bR$
with the situations that we will meet shortly over $\bC$.  The
manifold in the statement of the following proposition will often be a
manifold covering the one that we are dealing with, just as
$\bR/2\bZ$ covers $\bR/\bZ$ in the discussion above.

\begin{proposition}
\label{prop8.4}
Let $M$ be a compact connected Riemannian manifold, with its usual
ordinary metric, and let $u$ be a smooth function from $M$ to $\bR^n$
such that $u \cdot u = 1$.  Define the projection-valued function $p$
on $M$ by $p(m) = \<u(m),u(m)\>_0$.  Then
\[
L(p) = L(u).
\]
\end{proposition}

\begin{proof}
Let $m \in M$ and let $X$ be a tangent vector at $m$.  Let $D$ denote
``total derivative'', so that $D_X$ denotes differentiation at $m$ in
the direction of $X$.  Then
\[
D_Xp = \<D_Xu,u(m)\>_0 + \<u(m),D_Xu\>_0.
\]
Since $u \cdot u = 1$ we have $u(m) \cdot (D_Xu) = 0$.  We state the next
step as a lemma for later reference.  It is easy to prove by arguments
similar to those used in the proof of Proposition~\ref{prop8.2}.

\begin{lemma}
\label{lem8.5}
Let $\cH$ be a Hilbert space, over $\bR$ or $\bC$ and let 
$v,w \in \cH$ with $\|v\| = 1$ and $\<v, w\> = 0$.  Let
\[
T = \<v,w\>_0 + \<w,v\>_0.
\]
Then $\|T\| = \|w\|$.
\end{lemma}

On applying this lemma we see that $\|D_Xp\| = \|D_Xu\|$.  Since 
\[
\|(Dp)(m)\| = \sup\{\|D_Xp\|: \|X\| \le 1\}
\]
and similarly for $\|(Du)(m)\|$, we see that they are equal. 
Consequently
$\|Dp\|_{\i} = \|Du\|_{\i}$.  But standard arguments (see the
discussion early in Section \ref{sec11} after Corollary~\ref{cor14.2})
show that $L(p) = \|Dp\|_{\i}$ and similarly for $L(u)$.
\end{proof}

We remark that Proposition \ref{prop8.4} is false for functions from
$M$ to $\bC^n$ because the phase of $u$ can vary while leaving $p$
fixed. 

Actually, Proposition \ref{prop8.4} is true for path-length metric
spaces. We show this in Appendix A.

We now return to the M\"obius-strip bundle, and apply to it the
observations made above.  The consequence of the observations is that,
for our present purposes, it suffices to work with smooth $u$'s and
projections.  So let $u$ be, as earlier, a smooth function from $\bR$
to $\bR^n$ with $u \cdot u = 1$ and $u(r+1) = -u(r)$.
As $r$ goes from $0$ to $1$ the vector $u(r)$ traces a curve on the
unit sphere of $\bR^n$ from $u(0)$ to its antipodal point $-u(0)$.
The length of this curve is $\int_0^1 \|u'(r)\|dr$.  But the shortest
path from $u(0)$ to $-u(0)$ will be along one of the great-circle
geodesics, and it will have length $\pi$.  Since 
$\|p'(r)\| = \|u'(r)\|$ by Lemma \ref{lem8.5}, it follows that
$\int_0^1 \|p'(r)\|dr \ge \pi$.  By the mean-value theorem there must
be at least one point, $r_0$, where $\|p'(r_0)\| \ge \pi$.  Thus 
$L(p) \ge \pi$.  By the approximations discussed above, it follows
that for any projection $p \in M_n(A)$ such that $\Xi \cong pA^n$ we
have $L(p) \ge \pi$.  Also, we can achieve $L(p) = \pi$ by choosing
$u$ such that $u(r)$ moves along a great circle at speed $\pi$.

When we want to use Theorem~\ref{th6.2}, we see that it is best to
keep the size of our matrices as small as possible.  The simplest
choice is then $u(r) = (\cos(\pi r),\sin(\pi r))$.  The components of
this $u$ form a standard module frame for $\Xi$. We summarize what we
have found by:

\begin{proposition}
\label{prop8.6}
For any $p \in M_n(A)$ which represents the M\"obius-strip bundle we
have $L(p) \geq \pi$. For any $n \geq 2$ we can find such a $p$ with
$L(p) = \pi$. 
\end{proposition}

We remark that for any positive integer $k$ the space of $\bR$-valued
continuous functions $\xi$ which satisfy
\[
\xi(r-k)  = -\xi(r)
\]
is an $A$-module for pointwise multiplication, and it is an
entertaining and instructive exercise to show that these modules are
projective, to determine which are free, to find standard module
frames, etc.  In fact, the same is true for the modules consisting of
functions satisfying
\[
\xi(r-k) = +\xi(r).
\]

We now want to illustrate another aspect of our theory by examining
briefly what happens when one changes the metric on the circle.  Our
discussion will be at the qualitative level, but with more effort it
could be made quantitative.

Consider a smooth embedding of $M = \bT$ into $\bR^2$ with its
Euclidean metric, and assume that the image is approximately two 
far-away disjoint round circles connected by a very narrow ``tube''.
Give $M$ the metric coming from restricting the ordinary Euclidean
metric from $\bR^2$ to this embedding.  Let $p_0 \in M_2(C(M))$ be the
projection for a M\"obius-strip bundle such that $p_0$ is a constant
function on one of the almost-circles and the tube, so that the twist
takes place over the other almost-circle.  Let $p_1 \in M_2(C(M))$ be
the projection for a M\"obius-strip bundle such that $p_1$ is a
constant function on the tube and on the almost-circle where $p_0$ has
its twist.  Then one can show that for a suitable constant, $c$,
depending on the specific choice of embedding (especially the
narrowness of the tube), $p_0$ and $p_1$ can have been chosen to have
$L(p_j) < c$ for $j = 0,1$ but there is no homotopy $\{p_t\}$ of
projections from $p_0$ to $p_1$ such that $L(p_t) < c$ for all $t$.
Thus the set of projections $p \in M_2(C(M))$ which represent the
M\"obius-strip bundle and have $L(p) < c$ has more than one
path-component, and from our metric point of 
view the different path components can be viewed as representing
genuinely different vector bundles over $M$.  Indeed, $M$ with its
given metric can be made very close, for Hausdorff distance in
$\bR^2$, to the disjoint union of two circles (as seen by cutting the
very narrow tube), and $p_0$ and $p_1$ then correspond to vector
bundles on the disjoint union which are a M\"obius-strip bundle on one
circle and a trivial bundle on the other, but in different ways.  One
can make many variations of the above example, involving embedding $M$
as a greater number of almost-circles connected by narrow tubes.

One might object that the metrics involved in these examples are not
path-length metrics.  But one can use the same idea to smoothly embed
a $2$-sphere into $\bR^3$ as a collection of far away disjoint round
almost-spheres connected by narrow tubes.  Then instead of putting on
$M$ the restriction of the ordinary Euclidean metric on $\bR^3$, one
equips $M$ with the Riemannian metric from the embedding, and then the
ordinary metric from the Riemannian metric.  This is a path-length
metric.  Finally, one can consider projections corresponding to
putting on the various almost-spheres line bundles which on
corresponding actual spheres would have various Chern classes (as
discussed in Section~\ref{sec13}).


\section{Approximate M\"obius-strip bundles}
\label{sec9}

We now use the M\"obius-strip bundle to further illustrate our earlier
considerations, in a quantitative way.  The circle $\bT$ will now play
the role of the larger metric space $Z$ of our earlier discussion, and
so we denote the circle $\bT$ by $Z$ for the rest of this section.
Let $m$ be a large positive integer, and let 
$X = \{j/m: 0 \le j \le m-1\} \sim \bZ(1/m)/\bZ$ where the $j$'s are
integers.  We view $X$ as a subset of $Z$, and equip $X$ with the
metric from $Z$.  We now let $q_1$ denote the specific projection $p$
determined as in the previous section in terms of the standard frame
$(\cos(\pi r),\sin(\pi r))$.  Then we let $p_1$ denote the restriction
of $q_1$ to $X$.  Let us determine $L(p_1)$.  If $v$ and $w$ are two
unit-length vectors such that the angle between them is $\th$, then it
follows from Proposition~\ref{prop8.2} that the norm-distance between
the projections along these vectors is $|\sin \th|$.  Thus for 
$0 \le j,k < m$ we have 
\[
\|p_1(j/m) - p_1(k/m)\| = \sin(\pi|j-k|/m),
\]
and from this it is not hard to see that
\[
L(p_1) = \frac {\sin(\pi/m)}{(1/m)} = \pi \frac
{\sin(\pi/m)}{(\pi/m)}.
\]
Notice that this approaches $\pi$ as $m$ goes to $+\infty$, consistent with the
fact that $L(q_1) = \pi$ as seen in the previous section.

Since $X$ is finite, every vector bundle over $X$ is specified (up to
isomorphism) just by giving the dimension of the fiber vector-space
over each point.  Our projection $p_1$ represents the real vector
bundle whose fiber at each point has dimension $1$.  But this vector
bundle is equally well represented by the projection $p_0$ defined by
$p_0(t) = \begin{pmatrix} 1 & 0 \\ 0 & 0 \end{pmatrix}$ for each $t
\in X$.  Note that $L(p_0) = 0$.  Let $q_0$ denote the projection for
$Z$ defined by $q_0(r) = \begin{pmatrix} 1 & 0 \\ 0 & 0 \end{pmatrix}$
for all $r \in Z$, so that $p_0$ is the restriction of $q_0$ to $X$.
Then $q_0$ represents the rank-$1$ trivial vector bundle over $Z$, and
this bundle is not isomorphic to the M\"obius-strip bundle determined
by $q_1$.

The reason that this situation is possible, from our metric
point-of-view, is that there is no path, $p$, of projections from
$p_0$ to $p_1$ such that $L(p_t)$ is sufficiently small for all $t$.
This can be seen directly, by examining what such a path  must do at
various neighboring points of $X$.  But let us instead apply the
general considerations given in Theorem~\ref{th6.4}.  (We will not
expect this to give as sharp an estimate as a direct argument would
give.)  Set $\e = (2m)^{-1}$ and note that $X$ is $\e$-dense in $Z$
for this $\e$.  Let $p$ be a path of projections from $p_0$ to $p_1$,
and let $N$ be a constant such that $L(p_t) \le N$ for all $t$.  Note
that $\e L(q_0) = 0$, while $\e L(q_1) = \pi/2m < 1/2$ as soon as 
$m \ge 4$.  Now $M_2^s(\bR)$ is a $3$-dimensional vector space, and
from theorem~$1.1$a of \cite{KT} we find that the projection constant
of any $3$-dimensional real Banach space is no greater than $3$.  (The
techniques of section~$7$ of \cite{Rf76} can be used to obtain the
precise value of $\cP\cC(M_2^s(\bR))$.)  We will apply
Theorem~\ref{th6.4}, but by the observation just made we can replace
$\l_3$ there with $6$.  Suppose now that $m > 42N$ so that 
$\e 6N < 1/14$.  
We conclude from Theorem~\ref{th6.4} that there exists a continuous path
of projections from $q_0$ to $q_1$.  But we know that this is not
possible since $q_0$ and $q_1$ determine non-isomorphic $A$-modules.
Consequently we must have $N \ge m/42$.  Now $L(p_1) < \pi$ while
$L(p_0) = 0$.  Thus if $m \ge 4 \cdot 42$ so that $N > 4 > \pi$, we
see that the collection of projections in $M_2(A)$ which can be
connected to $p_1$ by paths $p$ of projections for which 
$L(p_t) \le 4$ for all $t$, does not include $p_0$.  Consequently the
projections in this collection can be viewed from our metric point of
view as giving approximate M\"obius-strip bundles on our finite set
$X$, which are not equivalent in our metric sense to the trivial
bundle of rank~$1$ on $X$.

Of course, similar considerations apply to other closed subsets of $Z$
which are $\e$-dense, and to other compact metric spaces $Y$ whose
Gromov--Hausdorff distance from $Z$ is less than $\e$ and for which a
corresponding metric on $Z  \dc  Y$ has been chosen.


\section{Lower bounds for $L(p)$ from Chern classes}
\label{sec10}

In this section we will indicate how Chern classes can sometimes be
used to obtain a lower bound on $L(p)$ for projections representing a
given vector bundle.  In Section~\ref{sec12} we will illustrate this
approach by considering vector bundles on a two-torus.  Our discussion
here is brief, and there is much more to be explored in this
direction.

For our purposes, and in particular for the two-torus, it is simplest
to work in the framework of Connes' 1980 paper \cite{Cn1}, which
initiated the subject of non-commutative differential geometry, and
which uses the Chern--Weil approach to Chern classes.  We briefly
sketch the setting.  We have a unital $C^*$-algebra $A$, together with
an action $\a$ of a connected Lie group $G$ by automorphisms of $A$.
(For our present purposes it will be quite sufficient for the reader
to have in mind just the case in which $A = C(G/H)$ where $H$ is some
cocompact closed subgroup of $G$, with $\a$ the evident action \cite{Rf77}.)  We
let $A^{\i}$ denote the dense $*$-subalgebra of smooth elements of $A$
with respect to $\a$.  Then $\a$ lifts to a homomorphism of the Lie
algebra, $\fG$, of $G$ (and its complexification) into the Lie algebra
$\mathrm{Der}(A^{\i})$ of derivations of $A^{\i}$ into itself.  We
denote this homomorphism again by $\a$.  We must also have a tracial
state, $\tau$, on $A$ which is invariant for the action 
$\a$.  For the case $A = C(G/H)$ this will just be a $G$-invariant 
probability 
measure on $G/H$ (unique if it exists).  

Let $\Xi$ be the smooth version of a projective $A$-module, that is, 
$\Xi$ is a projective $A^{\i}$-module.  On $\Xi$ there always exists a 
connection (i.e., covariant derivative), that is, a linear map 
$\nabla: \fG \to \mathrm{Lin}(\Xi)$ which satisfies the Leibniz
property
\[
\nabla_X(\xi a) = (\nabla_X(\xi))a + \xi(\a_X(a))
\]
for $X \in \fG$.  The curvature of $\nabla$ is the alternating 
$2$-form $\Th$ on $\fG$ defined by
\[
\Th(X,Y) = \nabla_X\nabla_Y - \nabla_Y\nabla_X - \nabla_{[X,Y]}.
\]
One finds that it has values in $\mathrm{End}_{A^{\i}}(\Xi)$.  Denote 
$\mathrm{End}_{A^{\i}}(\Xi)$ by $E$.

Let $\bw(\fG')$ denote the complexified exterior algebra over the dual 
vector space, $\fG'$, of $\fG$.  Then $\O^* = E \otimes \bw(\fG')$ has
in a natural way the structure of a differential graded algebra, and
we can view $\Th$ as an element of $\O^2$.  Thus 
$\Th \wedge \dots \wedge \Th$ 
($k$-times) can be viewed as an element of $\O^{2k}$, by using the
product 
in $E$.  The tracial state $\tau$ on $A$ induces in a natural way an 
unnormalized trace on $E$, which we denote by $\tau_E$.  It is 
characterized \cite{Rf33} by the property that
\[
\tau_E(\<\xi,\eta\>_E) = \tau(\<\eta,\xi\>_A)
\]
for $\xi,\eta \in \Xi$, where $\<\xi,\eta\>_E$ is the element 
(a ``rank-one'' operator) of $E$ defined by 
$\<\xi,\eta\>_E \z = \xi\<\eta,\z\>_A$.  Then 
$\tau_E (\Th \wedge \dots \wedge \Th)$, defined in the evident way, is
an 
element of $\bw^{2k}(\fG')$.  The main theorem \cite{Cn1} is that this 
element is a closed form, and that its cohomology class, $ch_k$,
depends 
only on $\Xi$.  (There are various choices of normalizing constants
that 
are used here.  We choose to use the constant $1$.)  If we pair $ch_k$ 
with a $2k$-homology class, we obtain a number.  This number may be
related 
to $L(p)$ when $p$ represents $\Xi$, but it is independent of the
choice of 
such (smooth) $p$.

Suppose now that we have a specific projection $p \in M_n(A^{\i})$
such that 
$\Xi = p(A^{\i})^n$.  On $(A^{\i})^n$ we have the evident flat
connection 
$\nabla^0$ given by $\nabla_X^0((a_j)) = (\a_X(a_j))$.  Then there is
a 
canonically associated connection, $\nabla$, on $\Xi$, defined by 
$\nabla_X(\xi) = p\nabla_X^0(\xi)$.  It is often called the Grassmann 
(or Levi--Civita) connection for $p$.  It is natural to use the
Grassmann 
connection in the setting sketched above.  The curvature of the
Grassmann 
connection is given \cite{Cn1} by
\[
\Th(X,Y) = p(\a_X(p)\a_Y(p) - a_Y(p)\a_X(p))p,
\]
where here $\a$ denotes the evident action of $G$ on $M_n(A)$.  
Clearly $E = \mathrm{End}_{A^\infty}(\Xi) = pM_n(A^\infty)p$.  Then $\tau_E$ is the 
restriction to $E$ of the canonical unnormalized trace $\tau$ on
$M_n(A)$ 
coming from $\tau$ on $A$.  In particular $\|\tau_E\| = \tau(p)$.  (We
remark 
that $\tau(p)$ is the $0$-th Chern class of $\Xi$.)  If we define $\o$
by
\[
\o(X,Y) = \tau_E(p(\a_X(p)\a_Y(p) - \a_Y(p)\a_X(p))),
\]
then $\o$ is a cocycle whose cohomology class, $ch_1$, is independent
of $p$ 
representing $\Xi$.  If we then pair $\o$ with a cycle in $\bw^2 \fG$,
then 
we obtain a number which is independent of $p$.  Now for any 
$X \wedge Y \in \bw^2\fG$ we have $d(X \wedge Y) = [X,Y]$.  
(See equation 3.1 of \cite{Ksz}.)  Thus $X \wedge Y$ is a cycle exactly 
if $[X,Y] = 0$.  Consequently, if $[X,Y] = 0$ then 
\[
c_{XY}(\Xi) = \tau_E(p(\a_X(p)\a_Y(p) - \a_Y(p)\a_X(p)))
\]
is a number independent of $p$ representing $\Xi$.

Suppose now that we have a norm,  $\nu$, on $\fG$, and that we define
a 
seminorm, $L$, on each $M_n(A)$ by
\[
L(a) = \sup\{\|\a_X(a)\|: \nu (X) \le 1\},
\]
for $\a$ extended to $M_n(A)$.  (See the next section.)  Then we find
that
\[
|c_{XY}(\Xi)| \le \|\tau_E\|2L(p)^2\nu (X)\nu (Y).
\]
Since $\|\tau_E\| = \tau(p)$, we obtain a lower bound for $L(p)$.
But by Theorem \ref{th3.5}, due to Hanfeng Li, this same lower
bound applies to any projection representing $\Xi$. Thus we obtain:

\begin{theorem}
\label{th10.1}
For every $p$ representing $\Xi$ we have
\[
(L(p))^2 \ge (2\tau(p))^{-1} \sup\{|c_{XY}(\Xi)|: [X,Y] = 0, \nu (X)
\le 1, \nu(Y) \le 1\}.
\]
\end{theorem}

To the extent that we have in hand cycles in $\bw^4\fG$, we can also
pair them with $\tau(\Th \wedge \Th)$ to obtain other lower-bounds for
$L(p)$, and similarly for higher dimensions.
More generally, for any ordinary compact Riemannian manifold, to the
extent that one has in hand specific even homology classes, one can
pair them with corresponding Chern classes of a given vector bundle
$\Xi$ to try to obtain lower bounds on $L(p)$ for smooth $p$'s
representing $\Xi$.  But consideration of flat bundles which are not
trivial shows that lower bounds from Chern classes may well not be
optimal bounds.

A possible related way of measuring the twisting of $\Xi$ would be just
by the size of the curvature $\Theta$ of various of its connections,
where by the size of $\Theta$ we mean 
$\sup\{\|\Theta(X, Y)\|: \nu(X), \nu(Y) \leq 1\}$ for a norm $\nu$
on $\fG$ as above. This measure of size is what is used in working
with ``almost flat bundles''. See section 3 of \cite{Skn} and the references
it contains. It would be interesting to investigate what could be
done in this direction. But again, there are non-trivial flat 
bundles, that have connections whose size measured this way is 0.


\section{Vector bundles on the two-torus}
\label{sec11}

We illustrate the considerations of the previous section by examining
the two torus.  One of our aims is to show that it can be feasible to
find projections $p$ with close-to-minimal $L(p)$.  We will let 
$G = \bT^2 = (\bR/\bZ)^2$ and $A = C(\bT^2)$, with the evident action
$\a$ of $G$ on $A$ by translation.  On $\fG = \bR^2$ we place the
standard inner product and corresponding norm.  This gives the
standard flat Riemannian metric on $\bT^2$, which in turn gives the
usual ordinary metric on $\bT^2$ given by
\[
d((r,s),(t,u)) = \min\{((r-t-m)^2 + (s-u-n)^2)^{1/2}: m,n \in \bZ\},
\]
where our notation on the left omits $\bZ^2$.  For
$a \in A^\infty$ the Lipschitz constant, $L(a)$, of $a$ can be conveniently 
calculated as
\[
L(a) = \sup\{\|\a_X(a)\|_{\i}: X \in \fG,\ \|X\| \le 1\},
\]
and similarly for $a \in M_n(A^\infty)$. See the proof of theorem 3.1 of 
\cite{Rf64} and lemma 3.1 of \cite{Rf70}.

The projective modules corresponding to the complex line bundles over 
$\bT^2$ can all be realized \cite{Rf40, Rf71} in the form
\[
\Xi_k = \{\xi \in C( \bT \x \bR , \bC): \xi(r,s+1) = e(kr)\xi(r,s)\},
\]
where $k \in \bZ$ and 
$e(t) = e^{2\pi it}$.  The action of elements of $A$ is by pointwise 
multiplication.

We need to find projections for these modules.  For this purpose it is 
clearer to work at first in somewhat greater generality.  Let $B$ be a 
unital $C^*$-algebra, and let 
\[
A = TB = \{a \in C( \bR , B): a(s+1) = a(s)\} \ .
\]
Our application will be to the case
in 
which $B = C(\bT)$ so that $A = C(\bT^2)$.  We use a variation of the 
familiar clutching construction to construct projective $A$-modules,
along 
the lines used in the proof of theorem~$8.4$ of \cite{Rf46}.  Let $w$
be a 
unitary element of $B$.  (In our application $w$ will be the function 
$e(kr)$.)  Set
\[
\Xi_w = \{\xi \in C(\bR , B): \xi(s+1) = w\xi(s)\} \ .
\]
Each $\Xi_w$ is a right 
$A$-module for pointwise multiplication, and has a natural Hermitian
metric 
given by $\<\xi,\eta\>_A(s) = \xi(s)^*\eta(s)$.  Furthermore, $\Xi_w$
is a 
projective $A$-module, for reasons that we now need to review.  
(See lemma~$8.8$ of \cite{Rf46}.)  We will apply
Proposition~\ref{prop7.2}, 
which will also give us a projection determining $\Xi_w$.  So we must
find a 
standard module frame for $\Xi_w$. Two elements of $\Xi_w$ suffice,
and we 
will denote them by $\eta_1$ and $\eta_2$.  We assume that $w$ is not 
connected to the identity through invertible elements, since otherwise 
$\Xi_w$ is easily seen to be a free module (lemma~$8.5$ of
\cite{Rf46}).  
It is reasonable to assume that $\eta_1(0) = 1_A$.  Since 
$\eta_1(s+1) = w\eta_1(s)$, we then see that $\eta_1$ must fail to be 
invertible at some point, since otherwise $w$ would be connected to
the 
identity.  In the absence of any further information about $w$ we will 
assume that $\eta_1$ actually takes value $0$ at some point $t_0$.  
Since $\|\eta_1(1)\| = \|\eta_1(0)\| = 1$, we see that to have a small 
Lipschitz norm it is best if $t_0 = 1/2$.  By the reconstruction
property, 
for any $\xi \in \Xi_w$ we must have
\[
\xi = \sum_{j=1}^2 \eta_j\<\eta_j,\xi\>_A = 
(\eta_1\eta_1^* + \eta_2\eta^*_2)\xi,
\]
so that we must have $\eta_1\eta_1^* + \eta_2\eta_2^* = 1_A$.  
It follows that $\eta_2(0) = 0$ and $\|\eta_2(1/2)\| = 1$.  Much as in 
\cite{Rf71} we will take $\eta_1$ and $\eta_2$ of the form
\[
\eta_1(s) = J_1(s)\cos (\pi s),\quad \eta_2(s) = J_2(s) \sin (\pi s)
\]
where
\begin{eqnarray*}
J_1(s) &= &(-w)^n \ \ \mbox{   for } \ n - 1/2 \ \le \ \ s \ < \  n +
1/2, \\
J_2(s) &= &(-w)^n \ \ \mbox{   for } \ n \ \le \ s \ < \ n + 1.
\end{eqnarray*}

Of course, $J_1$ and $J_2$ are discontinuous, but their
discontinuities 
are at the points where $\cos (\pi s)$ or $\sin (\pi s)$ takes value
$0$, 
and $\eta_1$ and $\eta_2$ are not only continuous, but are actually 
Lipschitz as functions on $\bR$.  For example, $\eta_2$ near $0$ is
given 
by $\eta_2(s) = \int_0^s h_2(t)dt$ where
\[
h_2(t) = \begin{cases}
\pi(-w^*)\cos(\pi t) &\mbox{for $-1 < t < 0$} \\
\pi\cos(\pi t) &\mbox{for $0 \le t < 1$ \ ,}
\end{cases}
\]
which is the derivative of $\eta_2$ where the derivative exists.  
If we let $\eta'_2$ denote this not-everywhere-defined derivative, we 
see that $L(\eta_2) = \|\eta'_2\|_{\i} = \pi$.  Similarly 
$L(\eta_1) = \|\eta'_1\|_{\i} = \pi$.  Since $J_j(s+1) = -w J_j(s)$
while 
$\cos(\pi s)$ and $\sin(\pi s)$ satisfy $g(s+1) = -g(s)$, we see that 
$\eta_j(s+1) = w\eta_j(s)$, so that $\eta_j \in \Xi_w$ for $j = 1,2$.  
Furthermore, we clearly have $\eta_1\eta_1^* + \eta_2\eta_2^* = 1_A$,
so 
that $\{\eta_1,\eta_2\}$ is a standard module frame for $\Xi_w$.  To
express 
the corresponding projection in $M_2(A)$ it is convenient to set 
$H = J_2^*J_1$, so that $H$ is the discontinuous periodic function
of period 1 taking 
value $1_A$ for $0 \le s < 1/2$ and value $-w$ for $1/2 \leq s < 1$.
Then, 
for example, $\<\eta_2,\eta_1\>_A(s) = H(s)\cos(\pi s)\sin(\pi s)$.
Thus
\[
p(s) = \begin{pmatrix}
\cos^2(\pi s)1_A &  \overline{H}(s)\cos(\pi s)\sin(\pi s) \\
H(s)\cos(\pi s)\sin(\pi s) &\sin^2(\pi s)1_A
\end{pmatrix} .
\]

We now apply all of this to the case in which 
$B = C(\bT)$ and $w(r) = e(kr)$ for a fixed $k$.  Then 
$A = C(\bT^2)$, and $p$ is a continuous projection-valued function on 
$\bT^2$ which is differentiable in $r$ and piecewise differentiable in
$s$.  Since $w(r) = e(kr)$, we see that $H$, as defined above 
but now for our special case, is given by
\[
H(s, r) = 1 \quad \mathrm{for} \quad 0 \leq s < 1/2 \quad 
\mathrm{and}  \quad -e(kr) \quad \mathrm{for} \quad 1/2 \leq s < 1 \ ,
\]
extended with period 1. Thus $H$ is continuously differentiable on $\bR^2$
whenever $s \notin \bZ /2$. From this we see that $p$ is continuously differentiable
on $\bR^2$ whenever $s \notin \bZ /2$. We now make the following
observation. Let $f$ be a continuous  and continuously piecewise-differentiable 
function, defined
on an interval $I$ in $\bR$, with values in a Banach space, such that at the points
of discontinuity of $f'$ the left-hand and right-hand limits exist. 
Then for any $s, t \in I$ we have
$f(t) - f(s) = \int_s^t f'(r)dr$ in the evident sense, given that $f'$ is not defined at a
finite number of points. Thus if there is a constant $K$ such
that $\|f'(r)\| \leq K$ when it is defined, then $L(f) \leq K$.
It is easily seen that for any two points $m, n \in \bR^2$ 
the restriction of $p$ to the straight line-segment joining those
two points satisfies the conditions on $f$ just stated. On applying
the above observation and using the periodicity of $p$, 
we see quickly that $L(p) \leq \| Dp \|_\infty$ where
$Dp$  denotes the total derivative of $p$, defined when $s \notin \bZ /2$.

We now express $p$ in terms of a unit vector field $u$. Most natural would
be to take $u = (\zeta_1, \zeta_2)'$, where ${}'$  denotes transpose.
But calculations are a bit simpler if we set 
$u(r,s) = (H(r, s)\cos(\pi s), \ \sin(\pi s))'$ and check that we still 
have $p = \langle u, u\rangle_0$, where, as in our
discussion 
of the M\"obius bundle, $\<u,u\>_0(r,s)$ denotes the rank-$1$
projection 
along the vector $u(r,s)$. Notice that $u$ is not even 
continuous where $s \in \bZ$. For any $X \in \bR^2$ let $D_X$ denote
the directional derivative along $X$ if it is defined. Then for 
$s \notin \bZ/2$ we have
\[
D_Xp = \<u,D_Xu\>_0 + \<D_Xu,u\>_0.
\]
Since $u \cdot u = 1$, we have $\mathrm{Re}(\<u,D_Xu\>_A) = 0$.  
Thus we need, here and later, the following small generalization of 
Lemma~\ref{lem8.5}, whose proof is obtained by using the 
techniques of the proof of Proposition~\ref{prop8.2}.

\begin{lemma}
\label{lem11.1}
Let $v$ and $w$ be vectors in a Hilbert space over $\bR$ or $\bC$, 
with $\|v\| = 1$ and $\mathrm{Re}(\<v,w\>_{\cH}) = 0$.  Let
\[
T = \<v,w\>_0 + \<w,v\>_0.
\]
Then $\|T\| = \|w - v\<v,w\>_{\cH}\| \le \|w\|$.
\end{lemma}

We thus see that
\[
\|D_Xp\|_{\i} = \|D_Xu - u\<u,D_Xu\>_{\cH}\|_{\i} \le \|D_Xu\|_{\i}.
\]
Now the total derivative, $Du$, of $u$, for $s \notin Z/2$,  is
\[
Du(r,s) = \begin{pmatrix}
\frac {\p H}{\p r} (r,s)\cos(\pi s) & -\pi H(r,s)\sin(\pi s) \\
0 & 
\pi \cos(\pi s)
\end{pmatrix} .
\]
This is complex, whereas we need the supremum of $\|D_Xu\|$ for $X$
real 
with $\|X\| = 1$.  By splitting $Du$ into its real and imaginary parts
it 
is not difficult to see that this supremum is $2\pi|k|$ for $k \ne 0$, 
obtained for $X = \begin{pmatrix} 1 \\ 0 \end{pmatrix}$ and for $s$ 
approaching $1$ from the left.  Thus $L(p) \le 2\pi|k|$ by 
Lemma \ref{lem11.1}.

Let us now obtain a lower bound for $L(p)$.  Straightforward
calculations 
show that, for the evident abbreviations,
\[
D_Xu - u\<u,D_Xu\>_{\cH} = \begin{pmatrix}
2\pi ikH\cos\sin^2 & -\pi H\sin \\
-2\pi ik \sin\cos^2 &\pi \cos
\end{pmatrix} \begin{pmatrix} X_1 \\ X_2 \end{pmatrix}
\]
for $X = \begin{pmatrix} X_1 \\ X_2 \end{pmatrix}$.  
At $s = 1/4$ and $r = 0$ this becomes 
\[
\frac {\pi}{\sqrt{2}} \begin{pmatrix}
ik & -1 \\
-ik & 1
\end{pmatrix} \begin{pmatrix}
X_1 \\ X_2 \end{pmatrix} = \frac {\pi}{\sqrt{2}} (ik X_1 - X_2) 
\begin{pmatrix} 1 \\ -1 \end{pmatrix},
\]
whose norm is $\pi(k^2X_1^2 + X_2^2)^{1/2}$.  For 
$\|X\| = 1$ this has maximum $\pi|k|$ for $k \ne 0$.  
Altogether we thus obtain
\[
\pi|k| \le L(p) \le \pi 2|k|.
\]
It is no surprise that $L(p)$ increases with $|k|$, since intuitively
the 
larger $|k|$ is the more rapidly $\Xi_k$ twists.  But further
investigation 
would be needed if one wanted to know whether $p$ is a projection for 
which $L(p)$ is minimal among all projections representing $\Xi_k$.  
In the next section we will see what information Chern classes can
give 
about this.

Let us now apply Theorem~\ref{th4.7}, using Notation~\ref{not4.6}, to
obtain 
the uniqueness of corresponding bundles on nearby spaces.  Let 
$X = \bT^2$ 
and let $\rho_X$ be the metric defined at the beginning of this
section.  
For each $k \in \bZ$ let $p^k$ be the projection for $\Xi_k$ defined
above.  
Suppose we have a metric space $(Y,\rho_Y)$ and a metric $\rho$ on 
$X  \dc  Y$ which restricts to $\rho_X$ and $\rho_Y$, such that 
$\mathrm{dist}_H^{\rho}(X,Y) < \e$.  Let $r$ be given with 
$r\e < 1/2$.  
For each $k$ such that $2\pi|k| < r$, so that $\e L(p^k) < 1/2$, 
let $\Pi_k$ denote the path component of $p^k$ in $\cP_2^r(X)$.  
Let $p'$ be another projection in $\Pi_k$, and let $s \ge r$ be such
that 
$\e s < 1/2$ and there is a path $p$ in $\Pi_k$ from $p^k$ to $p'$
such that 
for each $t$ there is a ${\bar q}_t \in \cP_2(Y)$ with 
$p_t \oplus {\bar q}_t \in \cP_2^s(X \dc Y)$.  Then 
according to  Theorem~\ref{th4.7} for any 
$q_0$ and $q_1 \in \cP_2^s(Y)$  such that $L(p^k \oplus q_0) < s$ and 
$L(p'  \oplus q_1) < s$, and for any $\d$ with $2\e s < \d < 1$, there
is 
a path ${\tilde p}$ from $p^k$ to $p'$ and a path $q$ from $q_0$ to
$q_1$ 
such that $L({\tilde p}_t \oplus q_t) < (1-\d)^{-1}s$ for each $t$.
In this 
controlled sense $q_0$ and $q_1$ are equivalent projections
corresponding 
to the equivalent projections $p^k$ and $p'$.

In the same way we can apply Theorem~\ref{th6.6} to obtain the
existence of 
vector bundles on nearby spaces, and in fact to lift homotopies.  
Let $X$ 
and $Y$ be as in the previous paragraph.  But now assume that 
$\e\l_2r < 1/12$.  Since $\l_2 \le 10/3$ by formula \eqref{eq5.3} and
the 
comments after Notation~\ref{not5.4}, it suffices to have 
$\e r < 1/40$.  
As in the previous paragraph let $\Pi_k$ be the path component of
$p^k$ in 
$\cP_2^r(X)$ for each $k$ such that $2\pi|k| < r$.  Set 
$s = \l_2r(1 - 12\e\l_2r)^{-1}$.  Then for any path $p$ in $\Pi_k$
there exists a path $q$ in $\cP^s(Y)$ such that 
$L^{\rho}(p_t \oplus q_t) < s$ for each $t$.


\section{Chern class estimates for the two-torus}
\label{sec12}

We now apply Theorem~\ref{th10.1} to obtain lower bounds for $L(p)$
for $p$'s representing the line bundles over $\bT^2$ discussed in the
previous section.  Thus, in the notation of Section \ref{sec10}, we
need to calculate $c_{XY}(\Xi_k)$, where now $X$ and $Y$ are in the
commutative Lie algebra $\bR^2$ of $\bT^2$.  To calculate this Chern
class we can choose any convenient smooth projection $p$ representing
$\Xi_k$.  We modify slightly our construction above for $p$ so as to
obtain a $p$ which is smooth.  We let $\var_1$ and $\var_2$ be
real-valued functions in $C^{\i}(\bR)$ which are close to 
$\cos(\pi s)$ and $\sin(\pi s)$ respectively, and satisfy both
$\var_j(s+1) = -\var_j(s)$ and $\var_1^2 + \var_2^2 = 1$, but are such
that $\var_1$ takes value $0$ in all of a small neighborhood of $1/2$,
while $\var_2$ takes value $0$ in all of a small neighborhood of $0$.
For $J_1$ and $J_2$ as above for our given $k$, we set 
$\eta_1(r,s) = J_1(r,s)\var_1(s)$ and 
$\eta_2(r,s) = J_2(r,s)\var_2(s)$.  Then $\eta_1$ and $\eta_2$ are
infinitely differentiable, and they form a standard module frame for
$\Xi_k$.  We let $p$ denote the corresponding projection, and much as
before we define $u$ by 
$u(r,s) = \begin{pmatrix} \eta_1(r,s) \\ \eta_2(r,s) \end{pmatrix}$.
Then, as before, we have $p = \<u,u\>_0$.

We want to express $c_{XY}(\Xi_k)$ in terms of $u$, and later, of
$\eta_1$ and $\eta_2$.  For $X \in \bR^2$ denote derivation in the $X$
direction by $\p_X$.  Then
\[
\p_X p = \<\p_Xu,u\>_0 + \<u,\p_Xu\>_0.
\]
But the $r$-component of $\p_Xu$ need not be a bounded function on
$\bR^2$, since $u$ satisfies $u(r,s+1) = e(kr)u(r,s)$.  We introduce
some notation to handle this.  Set $B = C^{\i}(\bR^2)$, the algebra of
smooth complex-valued, possibly unbounded, functions on $\bR$.  Thus
we can view (the smooth part of) $\Xi_k$ as a subspace of $B$.
Clearly $u \in B^2$, and $p \in M_2(B)$.  Actually much of our
calculation below works for any standard module frame for $\Xi_k$ whose
elements are smooth.  Thus we can assume just that $u \in B^n$ for
some $n$ (with entries in $\Xi_k$) 
so that $p \in M_n(B)$.  The main simplification which we
will be using is that we are dealing with a line bundle, so that the
$B$-valued (equivalently $A$-valued) inner product on $\Xi_k$ is given
just by pointwise multiplication.  If instead we were dealing with the
higher-rank (smooth) projective modules
\[
\Xi_{q,k} = \{\xi \in B: \xi(r,s+q) = e(kr)\xi(r,s), \ \xi(r+1,s) =
\xi(r,s)\}
\]
for a $q \geq 2$, then the inner product would involve a sum over
$\bZ/q\bZ$ in order to make its values periodic of period~$1$ in the
$s$ variable.  This would complicate matters.  Because we are in the
line-bundle situation, we can identify $E = \mathrm{End}_A(\Xi_k)$
with $A$ itself, and we have 
$\<\xi,\eta\>_E = \<\eta,\xi\>_A = \<\eta,\xi\>_B$ as functions.
(Note the reversal of order.)  Thus the key property of $u$ can be
rewritten as $\<u,u\>_B = 1$.  For $v,w \in B^n$ we let 
$\<v,w\>_0$ denote the corresponding ``rank-$1$'' operator defined by
$\<v,w\>_0x = v\<w,x\>_B$ for $x \in B^n$.  Then
$\<v,w\>_0 \in M_n(B)$, with $(\<v,w\>_0)_{jk} = \<w_k,v_j\>_B$.

For any $X \in \bR^2$ we will have the Leibniz rule
\[
\p_X(\<v,w\>_0) = \<\p_Xv,w\>_0 + \<v,\p_Xw\>_0,
\]
and similarly for $\<v,w\>_B$.
We now begin calculating, using strongly the line-bundle aspect
for $u$.  
Note first that
\[
p\<\p_Xu,u\>_0 = \<\<u,u\>_0\p_Xu,u\>_0 = \<u,\p_Xu\>_B \ p.
\]
It follows that
\[
p(\p_Xp) = p(\<\p_Xu,u\>_0 + \<u,\p_Xu\>_0) = 
\<u,\p_Xu\>_B \ p + \<u,\p_Xu\>_0.
\]
Now from $p^2 = p$ and the Leibniz rule we have $p(\p_Yp)p = 0$.  
Then, since $(\p_Yp)p$ is the adjoint of $p(\p_Yp)$, we have for 
$X,Y \in \bR^2$
\begin{eqnarray*}
p(\p_Xp)(\p_Yp)p &= &(\<u,\p_Xu\>_B \ p + \<u,\p_Xu\>_0)(\p_Yp)p \\
&= &\<u,\p_Xu\>_0(\p_Yp)p = \<u,\p_Xu\>_0(p\<\p_Yu,u\>_B +
\<\p_Yu,u\>_0) \\
&= &\<\<u,\p_Xu\>_0u,u\>_0\<\p_Yu,u\>_B + \<\<u,\p_Xu\>_0(\p_Yu),u\>_0
\\
&= &(\<\p_Xu,u\>_B\<\p_Yu,u\>_B + \<\p_Xu,\p_Yu\>_B)p.
\end{eqnarray*}
When we subtract from this the corresponding expression with $X$ and
$Y$ 
interchanged, we find that
\[
p(\p_X(p)\p_Y(p) - \p_Y(p)\p_X(p))p = 
2i\ \mathrm{Im}(\<\p_Xu, \ \p_Yu\>_B)p.
\]
The right-hand side is in $M_n(A)$ since the left-hand side is.  
The translation-invariant trace $\tau$ on $M_n(A)$ is given by taking
the 
ordinary trace of matrices followed by integration over the
fundamental 
domain $[0,1)^2 \sim \bT^2$.  Since $p$ is a rank-$1$ projection at
each 
point of $\bT^2$, we have $\tau(p) = 1$.  Thus when we apply the
formula 
for $c_{X,Y}$ given somewhat before Theorem \ref{th10.1} we find that
\[
c_{X,Y}(\Xi_k) = 2i \int_{\bT^2} \mathrm{Im}(\<\p_Xu,\p_Yu\>_B).
\]
(We are using the fact that our Lie algebra is commutative.)

For our specific $u$ defined in terms of $\eta_1$ and $\eta_2$ we then 
find that
\[
c_{X,Y}(\Xi_k) = 2i \int_{\bT^2}
\mathrm{Im}((\p_X\eta_1)^-(\p_Y\eta_1) + 
(\p_X\eta_2)^-(\p_Y\eta_2)).
\]
Let us take the particular case in which $X = \p/\p_r$ and 
$Y = \p/\p_s$, and 
denote these by $\p_1$ and $\p_2$ respectively.  Then for $j = 1,2$ we
have 
$\p_1\eta_j = (\p_1J_j)\var_j$ and $\p_2\eta_j = J_j(\p_2\var_j)$,
since 
$J_j$ is locally constant in $s$ and $\var_j$ vanishes in a
neighborhood of 
the discontinuities of $J_j$.  Examination of the definitions of $J_1$
and 
$J_2$ shows that
\[
\p_1J_1(r,s) = \begin{cases}
0 &\mbox{for $0 < s < 1/2$} \\
-2\pi ike(kr) &\mbox{for $1/2 < s < 1$,}
\end{cases}
\]
while $\p_1J_2(r,s) = 0$ for $0 < s < 1$.  Consequently 
$\p_1\eta_2 = 0$, 
and since $\p_1\eta_1 = (\p_1J_1)\var_1$ while 
$\p_2\eta_1 = J_1(\p_2\var_1)$, we find that
\begin{eqnarray*}
c_{\p_1,\p_2}(\Xi_k) &= &2i \int_{1/2}^1 
\mathrm{Im}((\p_1J_1)^-\var_1J_1(\p_2\var_1))ds \\ &=& 
4\pi ik \int_{1/2}^1 (1/2)\p_1(\var_1^2)ds 
= 2\pi ik(\var_1^2(1) - \var_1^2(1/2)) = 2\pi ik.
\end{eqnarray*}
We assume that $L$ is defined in terms of the standard inner product
on our 
Lie algebra $\bR^2$.  Since $\p_1$ and $\p_2$ come from elements of
our Lie 
algebra that commute and have norm $1$, we see from
Theorem~\ref{th10.1} 
that $L(p) \geq (\pi |k|)^{1/2}$. From Theorem \ref{th3.5},  due
to Hanfeng Li, we then obtain:

\begin{proposition}
\label{prop12.1}
For every projection $p$ representing $\Xi_k$ (of
any size) we have
\[
L(p) \ge (\pi|k|)^{1/2}.
\]
\end{proposition}

This estimate perhaps is not optimal, but it does show that as
$|k|$ 
goes to $+\i$ the lower bound for the $L(p)$'s goes to $+\i$, which is 
interesting, and consistent with what we found in the previous
section.  


\section{Projections for monopole and induced bundles}
\label{sec13}

In this section we begin to consider the complex line-bundles over the
$2$-sphere $S^2$, though most of our discussion will take place in
more general contexts that can be useful in treating other examples,
notably the coadjoint orbits of compact Lie groups considered in
\cite{Rf73}, as well as instantons \cite{LnS}.  See also section 2
of \cite{Rf77}. The term ``monopole
bundles'' traditionally refers to the non-trivial complex line-bundles
on $S^2$.  Our eventual aim is to understand what should be meant by
``monopole bundles'' on spaces close to $S^2$, in analogy with the
non-commutative case considered by physicists (for which see the
references at the beginning of the introduction).  From the action of
$SO(3)$ on $\bR^3$ we can view $S^2$ as $SO(3)/SO(2)$.  By means of
the adjoint representation of $SU(2)$ one obtains a $2$-sheeted
covering of $SO(3)$ by $SU(2)$, and through this covering we can also
view $S^2$ as $SU(2)/U(1)$, where $U(1)$ is the diagonal subgroup of
$SU(2)$.  Set $G = SU(2)$ and $H = U(1)$.  We can view $H$ as
$\bR/\bZ$, and then for each $n \in \bZ$ we can view the function 
$s \mapsto e(ns)$ as a character of $H$.  With this understanding, we
set 
\[
\Xi_n = \{\xi \in C(G): \xi(xs) = {\bar e}(ns)\xi(x) \mbox{ for } x
\in G,\ s \in H\}.
\]
With pointwise operations $\Xi_n$ is clearly a $C(G/H)$-module.  We
will see the well-known fact that it is projective, so corresponds to
a complex vector bundle, and in fact to a line-bundle. For $n \neq 0$
these are the monopole bundles. By definition their ``charge'' is
$|n|$.  Clearly $\Xi_n$ is carried into itself by the action of $G$ on
functions on $G$ by left translation, reflecting the fact that the
corresponding vector bundle is $G$-equivariant.

We seek projections for the $\Xi_n$'s.  The feature that we will use
to obtain projections is the well-known fact that the one-dimensional
representations of $H$ occur as subrepresentations of the restrictions
to $H$ of finite-dimensional unitary representations of $G$. Since $H
= U(1)$ is a maximal torus in $SU(2)$, the one-dimensional
representations which occur on restricting a representation of $G$ are
the weights of that representation.
We will treat the more general situation in which we have some compact
group $G$, a closed subgroup $H$, a unitary representation $V$ of $H$,
and a finite-dimensional unitary representation $(U,\cK)$ of $G$ whose
restriction to $H$ contains $V$ as a subrepresentation. Our approach
generalizes that given by Landi for the case of 2-spheres and
4-spheres in \cite{Lnd1, Lnd2}. (See also appendix A of \cite{LnS}.) 

The above set-up means that there is a subspace, $\cH$, of $\cK$ that
is carried into itself by the restriction of $U$ to $H$, and such that
this restricted representation of $H$ is equivalent to $V$.  From now
on we simply let $V$ denote this restricted representation.
In the next few paragraphs we work with continuous functions, but if
$G$ is a Lie group then everything has a version involving smooth
functions.  Set
\[
\Xi_V = \{\xi \in C(G,\cH): \xi(xs) = V_s^*(\xi(x)) \mbox{ for } x \in
G,\ s \in H\}.
\]
Clearly $\Xi_V$ is a module over $A = C(G/H)$ for pointwise
operations.  Also, $\Xi_V$ is clearly carried into itself by the
action of $G$ by left translation, i.e. $\Xi_V$ is $G$-equivariant.
(The completion of $\Xi_V$ for the inner product determined by a
$G$-invariant measure on $G/H$ is the Hilbert space for the unitary
representation of $G$ induced from the representation $V$ of $H$.)  On
$\Xi_V$ we define an $A$-valued inner product by 
\[
\<\xi,\eta\>_A(x) = \<\xi(x),\eta(x)\>_{\cK}.
\]
(We take the inner product on $\cK$, and so on $\cH$, to be linear in
the second variable.)

We want to show that $\Xi_V$ is a projective $A$-module, and find a
projection representing it.  Let $\O_{\cK} = C(G/H,\cK)$.  Then
$\O_{\cK}$ is a free $A$-module in the evident way.  For 
$\xi \in \Xi_V$ set $(\Phi\xi)(x) = U_x\xi(x)$ for $x \in G$, and
notice that $(\Phi\xi)(xs) = (\Phi\xi)(x)$ for $s \in H$ and 
$x \in G$, so that $\Phi\xi \in \O_{\cK}$.  It is clear that $\Phi$ is
an injective $A$-module homomorphism from $\Xi_V$ into $\O_{\cK}$.  We
show that the range of $\Phi$ is projective by exhibiting the
projection onto it from $\O_{\cK}$.  When applied to our earlier
$\Xi_n$ it will give a projection that we can then use in later
sections.  Let $P$ be the projection from $\cK$ onto $\cH$.  Note that
$U_sPU_s^* = P$ for $s \in H$ by the invariance of $\cH$.  Let $\cE$
denote the $C^*$-algebra $C(G/H,\cB(\cK))$, where $\cB(\cK)$ denotes
the algebra of operators on $\cK$.  In the evident way 
$\cE = \mathrm{End}_A(\O_{\cK})$.  Define $p$ on $G$ by 
\[
p(x) = U_xPU_x^*,
\]
and notice that $p(xs) = p(x)$ for $s \in H$ and $x \in G$, so that 
$p \in \cE$.  Clearly $p$ is a projection in $\cE$.  If 
$\xi \in \Xi_V$, then 
$p(x)(\Phi\xi)(x) = U_xPU_xU_x^*\xi(x) = (\Phi\xi)(x)$, so that
$\Phi\xi$ is in the range of $p$.  Suppose, conversely, that 
$F \in \O_{\cK}$ and that $F$ is in the range of $p$.  Set 
$\eta_F(x) = U_x^*F(x) = U_x^*p(x)F(x) = PU_x^*F(x)$.  Then the range
of $\eta_F$ is in $\cH$, and we see easily that 
$\eta_F(xs) = U_s^*\eta_F(x)$.  Thus $\eta_F \in \Xi_V$.  Furthermore,
$(\Phi\eta_F)(x) = F(x)$.  Thus $F$ is in the range of $\Phi$.  This
shows that the range of $p$ as a projection on $\O_{\cK}$ is exactly
the range of $\Phi$.  Thus this range, and so $\Xi_V$, are projective
$A$-modules.  Furthermore, $p$ is a projection for $\Xi_V$, so we have
attained our goal of finding a projection for $\Xi_V$.

To express $p$ as an element of $M_n(A)$ for some $n$ we need only
choose an orthonormal basis, $\{e_j\}_{j=1}^n$, for $\cK$, and view it
as a basis (so standard module frame) for $\O_{\cK}$, and express $p$
in terms of this basis.  Furthermore, if we define $\eta_j$ on $G$ by
$\eta_j(x) = PU_x^*e_j$, then it is easily seen that each $\eta_j$ is
in $\Xi_V$, and that $\{\eta_j\}$ is a standard module frame for
$\Xi_V$.  The basis also gives us an isomorphism of $\cE$ with
$M_n(A)$.

We remark that there will usually be many representations of $G$ whose
restriction to $H$ contains $V$ as a subrepresentation, and each of
these representations will give a projection for $\Xi_V$.  From the
statements of our earlier uniqueness and existence theorems we see
that it is probably best to take $n$ as small as possible, that is, to
take $\cK$ of as small dimension as possible.  Also, as we will see in
Section~\ref{sec16}, there may well be projections for $\Xi_V$ of even
smaller size that do not come from our construction above.

Finally, we remark that when $\cH$ is one-dimensional, $p$ is very
close to being a coherent state.  See the discussion in section~$3$ of
\cite{Rf73}.


\section{Metrics on homogeneous spaces}
\label{sec14}

Homogeneous spaces are quotient spaces, and the metrics that we will
use on them are quotients of metrics on the big space.  For technical
reasons we will need a strong understanding of the relation between a
quotient metric and the metric on the big space from which it comes,
both for ordinary metrics and for Riemannian metrics.  The definition
of (ordinary) quotient metrics is somewhat complicated.  A nice
exposition is given in \cite{Wvr2}, where three equivalent definitions
of quotient metrics are considered.  One of these definitions fits
especially well with our methods, and it goes as follows.  Let
$(Z,\rho)$ be a metric space, and let $\sim$ be an equivalence
relation on $Z$.  Let $A = C(Z)$, and define a closed subalgebra of
$A$ by 
\[
A/\sim \ \ = \ \{f \in C(Z): f(z) = f(w) \mbox{ if } z \sim w\}.
\]
(The functions in $C(Z)$ can be either real-valued or complex-valued.)
Let $L^{\rho}$ be the Lipschitz seminorm on $C(Z)$.  Define a
pseudometric, ${\tilde \rho}$, on $Z$ by
\[
{\tilde \rho}(z,w) = \sup\{|f(z)-f(w)|: f \in A/\sim \mbox{ and }
L^{\rho}(f) \le 1\}.
\]
Note that if $z \sim w$ then ${\tilde \rho}(z,w) = 0$.  Thus 
${\tilde \rho}$ drops to a pseudometric on the quotient space
$Z/\sim$.  But even if $Z$ is compact and all the equivalence classes
are closed (so that $Z/\sim$ with the quotient topology is Hausdorff
compact) ${\tilde \rho}$ may still fail to be a metric, and the
quotient metric space is obtained by identifying points which are at
distance $0$ for ${\tilde \rho}$.  However (see 1.16$_+$ of
\cite{Grm2}),

\begin{proposition}
\label{prop14.1}
Let $(Z,\rho)$ be a compact metric space, and let $G$ be a compact
group with an action $\a$ on $Z$, so that the quotient space $Z/\a$ is
compact and Hausdorff for the quotient topology.  If the action $\a$
is by isometries, then the corresponding pseudometric, 
${\tilde \rho}$, on $Z/\a$ is in fact a metric, giving the quotient
topology.  Let $\pi: Z \to Z/\a$ be the quotient map.  Then for any 
$f \in C(Z/\a)$ we have
\[
L^{\tilde \rho}(f) = L^{\rho}(f \circ \pi).
\]
\end{proposition}

\begin{proof}
Let $\cO_z$ be the $\a$-orbit of a point $z \in Z$, and define $f_z$
by $f_z(w) = \rho(w,\cO_z)$, for the evident meaning.  Then from the
fact that $\a$ is an action by isometries it is easily seen that $f_z$
is $\a$-invariant and so constant on orbits.  Furthermore it separates
$\cO_z$ from any other orbit, and $L^{\rho}(f_z) \le 1$.  From this it
follows easily that ${\tilde \rho}$ is a metric which gives the
quotient topology.  By composing with $\pi$ it is natural to view
$C(Z/\a)$ as a (closed $*$-) subalgebra of $C(Z)$, and then it
consists of all of the functions in $C(Z)$ which are constant on all
orbits.  For any $z \in Z$ and any $f \in C(Z/\a)$ we have 
$f(\pi(z)) = f(\cO_z)$, and from this it follows readily that
$L^{\tilde \rho}(f) = L^{\rho}(f \circ \pi)$.
\end{proof}

By applying linear functionals we quickly obtain:

\begin{corollary}
\label{cor14.2}
Let $(Z,\rho)$, $\a$, and $(X/\a,{\tilde \rho})$ be as above.  For any
Banach space $\cV$ and any function $F$ from $X/\a$ to $\cV$ we have
\[
L^{\tilde \rho}(F) = L^{\rho}(F \circ \pi).
\]
\end{corollary}

We now consider homogeneous spaces and Riemannian metrics.  Let $G$ be
a compact connected semisimple Lie group, with Lie algebra $\fG$.  On
$\fG$ we choose an Ad-invariant inner product, for example the
negative of the Killing form.  We denote this inner product by
$\<\cdot,\cdot\>_{\fG}$.  By translation this inner product defines a
bi-invariant Riemannian metric on $G$, with corresponding ordinary
metric $\rho$.

Let $\cV$ be a Banach space, over $\bR$ or $\bC$, and let $F$ be a
smooth function from $G$ into $\cV$.  Let $DF$ denote the total
derivative of $F$ for left translation, so that at each point 
$x \in G$ we can view $D_xF$ as a real-linear operator from $\fG$ into
$\cV$, defined by 
\[
(D_xF)(X) \ = \ (d/dt)|_{t=0} F(\exp(-tX)x).
\]
Simple arguments similar to those given in the proofs of theorem~$3.1$
of \cite{Rf64} and lemma~$3.1$ of \cite{Rf70} show that
\[
L^{\rho}(F) = \|DF\|_{\i} = \sup\{\|D_xF\|: x \in G\},
\]
where on $\fG$ we use the norm from its inner product.  (Our
Proposition~\ref{prop2.5} is also helpful here.) Similar
considerations appeared already early in Section \ref{sec11}.

Suppose now that $H$ is a closed subgroup of $G$, with Lie algebra
$\fH \subseteq \fG$.  Let $G/H$ be the corresponding homogeneous
space.  We can view functions in $C(G/H)$ as $H$-invariant functions
on $G$.  Now $G/H$ inherits a structure of differentiable manifold
\cite{Wrn} and we need to relate the tangent space at points of $G/H$
to those at points of $G$.  For this it is convenient to use an
explicit description of the projective module of smooth cross-sections
of the tangent bundle of $G/H$, along the lines given in \cite{Slb,
Frd, BIL, Rf77}.  This fits very well into the approach we have given
earlier concerning projective modules for vector bundles. Considerable
guidance in treating the tangent bundle of $G/H$ can also be found,
for example, in the proof of proposition~$3.16$ of \cite{CE} and its
surrounding discussion.  (A substantial part of our discussion can
easily be generalized to the more general setting of chapter~$3$ of
\cite{CE}, but that would take us beyond compact groups and spaces.)

We begin by letting $\fM$ denote the orthogonal complement to $\fH$ in
$\fG$.  Then $\fM$ is invariant under $\mathrm{Ad}_s$ for $s \in H$,
and $\fM$ can be identified with the tangent space at $eH \in G/H$.
At $xH$ each element of $\mathrm{Ad}_x(\fH)$ acts as the zero
derivation, and the tangent space at $xH$ can be identified with
$\mathrm{Ad}_x(\fM)$.  But one must be careful with this, since the
function $\mathrm{Ad}_{xs}(Y)$ for $Y \in \fM$ will not in general be
constant in $s$, hence the specific identification of
$\mathrm{Ad}_x(\fM)$ with the tangent space at $xH$ depends on the
choice of the coset representative $x$.  This can be handled by the
following description of the module $\cT(G/H)$ of smooth
cross-sections of the tangent bundle of $G/H$:

\begin{proposition}
\label{prop14.3}
The cross-section module $\cT(G/H)$ has a natural realization as 
\[
\{Z \in C^{\i}(G,\fG): Z(x) \in \mathrm{Ad}_x(\fM)
\mbox{  and  } Z(xs) = Z(x) \mbox{  for  } x \in G, s \in H\} .
\]
Given $F \in C^{\i}(G/H,\cV)$ for some Banach space $\cV$, the action
of $Z \in \cT(G/H)$ on $F$ is given by
\[
(ZF)(x) = (d/dt)|_{t=0} F(\exp(-tZ(x))x).
\]
Equivalently, the cross-section module could be given by
\[
\{W \in C^{\i}(G,\fM): W(xs) = \mathrm{Ad}_s^{-1}(W(x))
\mbox{  for all  } x \in G, s \in H\} ,
\]
with the action of $W$ on $F$ given by
\[
(WF)(x) = (d/dt)|_{t=0} F(x\exp(-tZ(x))).
\]
\end{proposition}

\begin{proof}
It is clear that $\cT(G/H)$ is a $C_{\bR}^{\i}(G/H)$-module for
pointwise operations.  Let $Q$ denote the orthogonal projection of
$\fG$ onto $\fM$.  Then $Q \mathrm{Ad}_s = \mathrm{Ad}_sQ$ for each 
$s \in H$.  Each $X \in \fG$ determines an element ${\tilde X}$ of
$\cT(G/H)$, defined by
\[
{\tilde X}(x) = \mathrm{Ad}_x(Q(\mathrm{Ad}_{x^{-1}} X))
\]
for $x \in G$.  Note that there may well be points $x$ where
$\mathrm{Ad}_{x^{-1}}X \in \fH$, so that ${\tilde X}(x) = 0$.  (We are
feeling the effect here of the fact that $G/H$ may well not be
parallelizable.)  However, we see that for any 
$X \in \mathrm{Ad}_x(\fM)$ we have ${\tilde X}(x) = X$, so that the
evaluations at $x$ of elements of $\cT(G/H)$ fill the whole tangent
space at $xH$.  From this and the fact that $\cT(G/H)$ is a
$C_{\bR}^{\i}(G/H)$-module it follows that $\cT(G/H)$ does represent
the full space of smooth cross-sections of the tangent bundle for
$G/H$.

When we use the fact that
\[
\exp((Z(y))x = x\exp(\mathrm{Ad}_x^{-1}(Z(y))
\]
we obtain the second description of the cross-section module.
\end{proof}

We will use the first description of $\cT(G/H)$ given by
the above proposition.
We define a $C_{\bR}^{\i}(G/H)$-valued inner product on $\cT(G/H)$ by
\[
\<Z,W\>_{G/H}(x) = \<Z(x),W(x)\>_{\fG}.
\]
Calculations very similar to those in the previous section show that
if $\{X_j\}$ is an orthonormal basis for $\fG$, then $\{{\tilde
X}_j\}$ is a standard module frame for $\cT(G/H)$.  (Here, of course,
we are working over $\bR$.)  Furthermore, the
$C_{\bR}^{\i}(G/H)$-valued inner product gives the Riemannian metric
on $G/H$ which is induced from that on $G$ as discussed in chapter~$3$
of \cite{CE}.

The Riemannian metric on $G/H$ will define an ordinary metric 
${\hat \rho}$ on $G/H$.  Given a smooth function $F$ on $G/H$ with
values in some Banach space, its Lipschitz constant, 
$L^{\hat \rho}(F)$, for the ordinary metric ${\hat \rho}$ on $G/H$ is
defined.  But, for the kinds of reasons discussed above for $\rho$, we
also have
\[
L^{\hat \rho}(F) = \sup\{\|ZF\|_{\i}: Z \in \cT(G/H) \mbox{ and }
\<Z,Z\>_{G/H} \le 1\}.
\]
However, we can also view $F$ as a function on $G$, and calculate its
Lipschitz constant $L^{\rho}(F)$ there.  As seen above, we have
\[
L^{\rho}(F) = \sup\{\|XF\|_{\i}: X \in \fG,\ \|X\| \le 1\},
\]
where the elements of $\fG$ are viewed as right-invariant vector
fields.  The relation that we need for later calculations is:

\begin{proposition}
\label{prop14.4}
Let $F \in C^{\i}(G/H,\cV)$ for some Banach space $\cV$, and view $F$
also as a function on $G$.  Then
\[
L^{\hat \rho}(F) = \sup\{\|XF\|_{\i}: X \in \fG,\ \|X\| \le 1\}.
\]
Consequently $L^{\hat \rho}(F) = L^{\rho}(F)$, and ${\hat \rho}$ is
the quotient metric from $\rho$.
\end{proposition}

\begin{proof}
Let $Z \in \cT(G/H)$ with $\<Z,Z\>_{G/H} \le 1$, so that 
$\|Z_x\| \le 1$ for all $x \in G$.  From Proposition~\ref{prop14.3} we
see that $\|(ZF)(x)\| = \|Z_xF\|$, and from this it is evident that
$L^{\hat \rho}(F) \le L^{\rho}(F)$.  On the other hand, let $x \in G$
be given, and let $X \in \mathrm{Ad}_x(\fM)$ with $\|X\| \le 1$.
Define ${\tilde X} \in \cT(G/H)$ as done earlier, and recall that
${\tilde X}(x) = X$.  It follows that 
$\|{\tilde X}F\|_{\i} \ge \|(XF)(x)\|$.  Notice from the definition of
${\tilde X}$ that $\<{\tilde X},{\tilde X}\>_{G/H} \le 1$.  It now
follows that $L^{\rho}(F) \le L^{\hat \rho}(F)$, so that they are
equal.  If we apply this for $F \in C^{\i}(G/H)$, we see that 
${\hat \rho}$ is the quotient of $\rho$.
\end{proof}


\section{Lipschitz constants of induced bundles}
\label{sec15}

We now combine the results of Sections~\ref{sec13} and \ref{sec14} to
obtain a formula for the Lipschitz constants of the projections
obtained in Section~\ref{sec13} when our setting is compact semisimple
Lie groups as in Section~\ref{sec14}.  Thus we let $G$, $\fG$, $H$,
etc., be as in Section~\ref{sec14}, and we now assume that we have a
unitary representation $(U,\cK)$ of $G$ and an $H$-invariant subspace
$\cH$ of $\cK$, with the restriction of $U$ to $H$ and $\cH$ denoted
by $V$.  As before we let $p(x) = U_xPU_x^*$, where $P$ is the
orthogonal projection of $\cK$ onto $\cH$.  Now finite-dimensional
representations of a Lie group are smooth, and so $p$ is a smooth
function.  We let $Dp$ denote its total derivative.  From the previous
section we see that $L(p) = \|Dp\|_{\i}$.  Now a simple calculation
shows that for any $x \in G$ and $X \in \fG$ we have
\[
(Dp)_x(X) = [U_xPU_x^*, \ U_X]
\]
where $U_X$ refers to the infinitesimal version of $U$ as a
representation of $\fG$ on $ \cK$.  Then
\[
\|(Dp)_x(X)\| = \|U_x[P,U_x^*U_XU_x]U_x^*\| 
= \|[P,U_{\mathrm{Ad}_{x^{-1}}(X)}]\|.
\]
Since we have chosen the inner product on $\fG$ so that it is
preserved by Ad, it follows that
\[
\|(Dp)_x\| = \sup\{\|[P,U_X]\|: \|X\|_{\fG} \le 1\}.
\]
Since the right-hand side is independent of $x$, we see that it also
equals $\|Dp\|_{\i}$.  But each $U_X$ is skew adjoint, and so
\[
[P,U_X] = PU_X - U_XP = PU_X + (PU_X)^* = 2\mathrm{Re}(PU_X).
\]
When we combine this with Proposition~\ref{prop14.4} we obtain:

\begin{proposition}
\label{prop15.1}
With notation as above,
\[
L_{G/H}(p) = \|Dp\|_{\i} = 2\sup\{\|\mathrm{Re}(PU_X)\|: X \in \fG,\
\|X\| \le 1\},
\]
where $L_{G/H}$ denotes the Lipschitz seminorm for the metric ${\hat
\rho}$ on $G/H$.
\end{proposition}

Suppose now that our representation $V$ of $H$ is one-dimensional, so
that $\Xi_V$ is the space of cross-sections for a line bundle.  Let
$v_0$ be a unit vector in the Hilbert subspace $\cH \subset \cK$ for
$V$, so that $P = \<v_0,v_0\>_0$, the rank-$1$ operator determined by
$v_0$.  Notice that since $U_X$ is a skew-adjoint operator,
$\mathrm{Re}(\<v_0,U_Xv_0\>_{\cH}) = 0$.  Then on combining
Proposition~\ref{prop15.1} with Lemma~\ref{lem11.1} we obtain:

\begin{proposition}
\label{prop15.2}
With notation as at the beginning of this section, but with $\cH$ of
dimension~$1$, spanned by the unit-vector $v_0$, we have
\[
L_{G/H}(p) = \sup\{\|(I-P)U_Xv_0\|: X \in \fG,\ \|X\| \le 1\}.
\]
\end{proposition}

For use in the next section we also need to consider tensor powers of
one-dimensional representations.  For $(U,\cK)$, $v_0$, $(V,\cH)$, $P$
and $p$ as above and for a given $n$ let $\cK^{\otimes n}$ be the $n$-fold full tensor
power of $\cK$, with $U^{\otimes n}$ the corresponding representation
of $G$.  Let $v_0^n$ be the $n$-fold tensor power of $v_0$, and let
$\cK_n$ be the $G$-invariant subspace of $\cK^{\otimes n}$ generated
by $v_0^n$.  Let $P_n$ be the projection of $\cK_n$ onto the
one-dimensional subspace, $\cH_n$, spanned by $v_0^n$, which is
$H$-invariant.  Let $U^{(n)}$ be the restriction of $U^{\otimes n}$ to
$\cK_n$, and let $V^n$ be the restriction of $U^{\otimes n}$ to a
representation of $H$ on $\cH_n$.  Finally, let 
$p_n \in C(G/H,B(\cK_n))$ be defined by 
$p_n(x) = U_x^{(n)}P_n(U_x^{(n)})^*$, so that we are in a version of
the setting discussed above.  Then
\[
L(p_n) = \sup\{\|(I-P_n)U_X^{(n)}v_0^n\|: X \in \fG,\ \|X\| \le 1\}.
\]
But
\[
U_X^{(n)}v_0^n = (U_Xv_0) \otimes v_0 \dots \otimes v_0 \quad + \quad
v_0 \otimes (U_Xv_0) \otimes v_0 \dots v_0 \quad + \quad \dots
\]
so that
\[
\<v_0^n,U_X^{(n)}v_0^n\> \ = \ n\<v_0,U_Xv_0\>.
\]
Let $w_X \ = \ U_Xv_0 - v_0\<v_0,U_Xv_0\>_{\cH}$.  Then
\begin{eqnarray*}
U_X^{(n)}v_0^n &- &v_0^n\<v_0^n,U_X^{(n)}v_0^n\> \ = \ U_X^{(n)}v_0^n -
nv_0^n\<v_0,U_Xv_0\>_{\cH} \\
&= &w_X \otimes v_0 \otimes \dots \otimes v_0 \ \ + \ \ v_0
\otimes w_X \otimes v_0 \otimes \dots \otimes v_0 \ \ + \ \ \dots
\end{eqnarray*}
Thus, since $w_X \bot v_0$, we find that
\[
\|U_X^{(n)}v_0^n \ - \ v_0^n\<v_0^n,U_X^{(n)}v_0^n\>\|^2 \ = \ n\<w_X,w_X\>.
\]
Combining this with Lemma~\ref{lem11.1}, we obtain:

\begin{proposition}
\label{prop15.3}
With notation as above,
\[
L(p_n) = n^{1/2}L(p).
\]
\end{proposition}


\section{The sphere, $SU(2)$, and monopole and induced bundles}
\label{sec16}

To treat the two-sphere $S^2$ we view it as $G/H$ where $G = SU(2)$
and $H$ is the diagonal subgroup 
$\begin{pmatrix} e(t) & 0 \\ 0 & {\bar e}(t) \end{pmatrix} =
\begin{pmatrix} z & 0 \\ 0 & {\bar z} \end{pmatrix}$, \  as discussed
at the beginning of Section~\ref{sec13}.  For each $n \in \bZ$ we
consider the representation $t \mapsto e(nt)$ of $H \cong \bR/\bZ$,
and the corresponding induced $A$-module $\Xi_n$, where $A = C(S^2)$.
We begin by considering the case $n = 1$.  The corresponding
representation occurs in the restriction to $H$ of the standard
representation of $SU(2)$ on $\bC^2$.  As $v_0$ we can take the
lowest-weight vector $\begin{pmatrix} 0 \\ 1 \end{pmatrix}$.  Then the
projection onto its span is 
$P = \begin{pmatrix} 0 & 0 \\ 0 & 1 \end{pmatrix}$.  For the
corresponding projection $p_1$ in $M_2(A)$ we have, by
Proposition~\ref{prop15.2},
\[
L(p_1) = \sup\{\|(I-P)U_Xv_0\|: X \in \fG,\ \|X\| \le 1\}.
\]
A general element of $\fG = su(2)$ will be of the form 
$X = \begin{pmatrix} ir & -{\bar w} \\ w & - ir \end{pmatrix}$ where
$r \in \bR$ and $w \in \bC$.  Then 
$\|(I-P)U_Xv_0\| = 
\left\| \begin{pmatrix} -{\bar w} \\ 0 \end{pmatrix}\right\| = |w|$.
We choose our normalization of the Ad-invariant inner product on
$su(2)$ to be 
$\<X,Y\>_{\fG} = 2^{-1} \mathrm{tr} (XY^*)$, so that for $X$ as above
we have $\|X\| = (r^2 + |w|^2)^{1/2}$.  For this choice we see that
$L(p_1) \le 1$.  But one choice of $X$ is 
$\begin{pmatrix} 0 & -1 \\ 1 & 0 \end{pmatrix}$, and from this choice
we see that $L(p_1) = 1$.  If instead we had chosen as $v_0$ the
highest-weight vector $\begin{pmatrix} 1 \\ 0 \end{pmatrix}$, we would
be dealing with the case $n = -1$, and so $\Xi_{-1}$.  All of the
above discussion applies with almost no change to this case too.  We
have thus obtained:

\begin{proposition}
\label{prop16.1}
Let $p$ be the projection in $M_2(C^{\i}(S^2))$ constructed above for
either of the monopole bundles $\Xi_1$ and $\Xi_{-1}$.  For the
Riemannian metric on $S^2$ coming from the choice of normalization of
the inner product on $su(2)$ made above, we have $L(p) = 1$.
\end{proposition}

Let us now apply the results of Sections~\ref{sec4} and \ref{sec6}.
We could do this in the same way as we did near the end of
Section~\ref{sec11} for the two torus.  But for the sake of variety,
let us take a slightly different approach, involving specific choices
of numbers, and more focused on our monopole projection above while
less focused on homotopies.  So let $p_*$ denote our projection above
for $\Xi_1$ or $\Xi_{-1}$.  For the moment let us denote $S^2$ by $X$
and let $\rho_X$ be its round metric corresponding to our choice above
of an Ad-invariant inner product on $su(2)$.

Choose $\e < 1/60$, so that $12\e\l_2L(p_*) < 2/3$, since 
$\l_2 \leq 10/3$ by equation \ref{eq5.3}.  Let $(Y,\rho_Y)$ be another
compact metric space, and suppose that we have a metric $\rho$ on 
$X \dc Y$ that restricts to $\rho_X$ and $\rho_Y$, and for which
$\mathrm{dist}_H^{\rho}(X,Y) < \e$.  Then according to
Corollary~\ref{cor6.3} there is a projection $q_1 \in M_2(C(Y))$ such
that
\[
L(p_* \oplus q_1) < \l_2L(p_*)(1 - (2/3))^{-1} \le (10/3)(3/1) =10.
\]
Then $\e L(p_* \oplus q_1) < 1/6 < 1/2$, and so from
Theorem~\ref{th4.5} we find that if $q_0$ is any other projection in
$M_2(C(Y))$ such that $L(p_* \oplus q_0) < 10$ then there is a path
$q$ of projections in $M_2(C(Y))$ going from $q_0$ to $q_1$ such that
\[
L(p_* \oplus q_t) < (1 - 20/60)^{-1}10 = 15.
\]
In particular $q_0$ and $q_1$ (and all the $q_t$'s) will determine
isomorphic vector bundles on $Y$, and so it is reasonable to consider
the corresponding isomorphism class of bundles on $Y$ to be the
``monopole'' bundle on $Y$ for the given proximity of $Y$ to $X$,
although the statement in terms of projections is more precise.

We now consider $\Xi_n$ for $|n| \ge 2$.  We carry out the discussion
for $n \ge 2$, but a very parallel discussion works for $n \le -2$.
We seek to apply Proposition~\ref{prop15.3}.  Let 
$v_0 = \begin{pmatrix} 0 \\ 1 \end{pmatrix} \in \bC^2 = \cH$, and let
$U$ be the standard representation of $SU(2)$ on $\cK = \bC^2$ as
above.  Form $U^{\otimes n}$ and $\cK^{\otimes n}$ as done before
Proposition~\ref{prop15.3}, and then form $v_0^n$, $\cH_n$ and
$U^{(n)}$.  Because $v_0$ is a lowest weight vector, it is not
difficult to see that $U^{(n)}$ on $\cK_n$ is the $(n+1)$-dimensional
irreducible representation of $SU(2)$ and that $v_0^n$ is a
lowest-weight vector for it.  (See for example proposition~VII.2 of
\cite{Srr} for a more general context.)  We clearly have
$U_s^{(n)}v_0^n = {\bar e}(ns)v_0^n$ for $s \in H$, so the restriction
of $U^{(n)}$ to $H$ and to the one-dimensional subspace, $\cH_n$,
spanned by $v_0^n$ is the representation of $H$ which defines $\Xi_n$.
Thus if, as before, we let $P_n$ denote the projection from $\cK_n$
onto $\cH_n$, and if we set $p_n(x) = U_x^{(n)}P_n(U_x^{(n)})^*$, then
$p_n$ is a projection in $C(G/H, \cB(\cK_n))$ which represents
$\Xi_n$.  From Proposition~\ref{prop15.3}, applied also for $n \le -2$
(and also from Proposition~\ref{prop16.1} for $n = \pm 1$, and for
$p_0$ a constant projection for the free rank-$1$ projective module
$\Xi_0$ over $G/H$), we obtain:

\begin{theorem}
\label{th16.2}
For any $n \in Z$ and for the projection $p_n$ defined above for
$\Xi_n$ we have
\[
L(p_n) = |n|^{1/2}.
\]
\end{theorem}

We notice however that since $\cK_n$ has dimension $|n| + 1$, the
standard module frame for $\Xi_n$ corresponding to $p_n$ and a basis
for $\cK_n$ will have $|n| + 1$ elements.  This corresponds to the
fact that in the present setting $\Xi_n$ is embedded in a free module
of rank $|n| + 1$.  But any complex line-bundle on a compact space of
dimension~$2$ can be obtained as a summand of a rank-$2$ trivial
bundle.  (See, for example, proposition~$1.1$ and theorem~$1.5$ of
chapter~$8$ of \cite{Hsm}.)  So we should be able to obtain $\Xi_n$ as
a summand of a free module of rank~$2$, with a frame consisting of
just two elements and a corresponding projection in $M_2(A)$.  We can
indeed do this, as follows.

Fix $n$, and let $\z_1$ and $\z_2$ be the functions on $SU(2)$ defined
by $\z_1(x) = {\bar z}_1^n$ and $\z_2(x) = {\bar z}_2^n$ respectively
when $n > 0$, and by their complex conjugates when $n < 0$, where much
as before 
$x = 
\begin{pmatrix} z_1 & -{\bar z}_2 \\ z_2 & {\bar z}_1 \end{pmatrix}$
with $|z_1|^2 + |z_2|^2 = 1$.  It is clear that $\z_1$ and $\z_2$ are
in $\Xi_n$.  From now on we assume that $n > 0$, but entirely parallel
considerations apply for $n < 0$.  Let $h \in A$ be defined by 
$h(x) = (|z_1|^{2n} + |z_2|^{2n})^{-1/2}$,  and set $\eta_j = \z_jh$
for $j = 1,2$, so that 
$\<\eta_1,\eta_1\>_A + \<\eta_2,\eta_2\>_A = 1$.  Since $\Xi_n$
corresponds to a line bundle, $\mathrm{End}_A(\Xi_n) = A$, and so this
last relation is basically the reconstruction formula of
Definition~\ref{def7.1}.  Thus $\{\eta_1,\eta_2\}$ is a standard
module frame for $\Xi_n$, and this frame provides a projection, $p$,
in $M_2(A)$ for $\Xi_n$.  Specifically, 
$p_{jk} = \<\eta_j,\eta_k\>_A = {\bar \eta}_j\eta_k$,
so that 
\[
p_{jk}(x) = z_j^n{\bar z}_k^n/(|z_1|^{2n} + |z_2|^{2n})
\]
for $j,k = 1,2$. If we define $u$ on $SU(2)$ by 
$u(x) = h(x) \begin{pmatrix} z_1^n \\ z_2^n \end{pmatrix}$, then 
$p(x) = \<u(x),u(x)\>_0$, the rank-$1$ operator on $\cK = \bC^2$
determined by the unit-vector $u(x)$.  (Recall that we take the inner
product on $\cK$ to be linear in the second variable.)  Note that $u$
is not constant on cosets of $H$, but that $p$ is.

We can now try to calculate the Lipschitz constant of $p$, first as a
function on $G$ and then on $G/H$ by using Proposition~\ref{prop14.4}
much as before.  It is clear that $p$ and $u$ are smooth.  Thus we can
use derivatives to calculate $L(p)$.  The advantage of working on $G$
is that while $u$ is not a function on $G/H$, we can make use of its
derivatives on $G$ to calculate those of $p$.

Suppose now that $u$ is any smooth function from $G$ into a Hilbert
space $\cK$ such that $\<u(x),u(x)\>_{\cK} = 1$ for all $x$, and
define $p$ by $p(x) = \<u(x),u(x)\>_0$.  Then for any $X \in \fG$ and
$x \in G$ we have, much as earlier,
\[
(Xp)(x) = \<(Xu)(x),u(x)\>_0 + \<u(x),(Xu)(x)\>_0.
\]
Because each $u(x)$ is a unit-vector we see that
$\mathrm{Re}(\<(Xu)(x),u(x)\>_{\cK}) = 0$.  For any fixed $x \in G$ we
can then apply Lemma~\ref{lem11.1} to conclude that
\[
\|(Xp)(x)\| = \|(Xu)(x) - u(x)\<u(x),(Xu)(x)\>\| = \|(I_{\cK} -
p(x))((Xu)(x))\|.
\]
Consequently, much as just after Lemma~\ref{lem11.1}, we obtain:

\begin{proposition}
\label{prop16.3}
Let $u$ be a smooth function from $G$ to a Hilbert space $\cK$ such
that $\<u(x),u(x)\>_{\cK} = 1$ for all $x \in G$.  Let $p$ be defined
by $p(x) = \<u(x),u(x)\>_0$.  Then
\[
L(p) = \sup\{\|(I_{\cK} - p(x))((Xu)(x))\|: x \in G,\ \|X\| \le 1\}.
\]
\end{proposition}

When the dimension of $\cK$ is $2$, the projection $I_{\cK} - p(x)$
will be of rank~$1$, and so is given by a unit vector, $v(x)$,
orthogonal to $u(x)$.  Then
\begin{eqnarray*}
&\|(I_{\cK}-p(x))((Xu)(x))\| = \|\<v(x),v(x)\>_0((Xu)(x))\| \\
& \quad \quad=\|v(x)\<v(x),(Xu)(x)\>_{\cK}\| = |\<v(x),(Xu)(x)\>_{\cK}|.
\end{eqnarray*}
For $u(x) = h(x)({\bar \z}_1(x),{\bar \z}_2(x))'$ as before, where the
prime denotes transpose, we can choose $v(x)$ to be defined by 
$v(x) = h(x)(\z_2(x),-\z_1(x))'$.  Now
\[
(Xu)(x) = (Xh)(x)u(x)(h(x))^{-1} \ + \ h(x)((X{\bar \z}_1)(x),(X{\bar \z}_2)(x))',
\]
so that, since $v(x) \bot u(x)$,
\begin{eqnarray*}
\<v(x),(Xu)(x)\>_{\cK} &= &\<v(x),h(x)((X{\bar \z}_1)(x),(X{\bar
\z}_2)(x))'\>_{\cK} \\
&= &h^2(x)(\z_2(x)(X\z_1)(x) - \z_1(x)(X\z_2)(x)).
\end{eqnarray*}
Since $\z_1$ and $\z_2$ are defined for all $(z_1,z_2) \in \bC^2$, and
as functions on $G$ they are determined by the first column of $x$, we
can calculate $(X\z_j)(x)$ by means of the straight-line path $I - tX$
instead of the path $\exp(-tX)$, for $t \in \bR$.  For 
$X = \begin{pmatrix} ir & -{\bar w} \\ w & -ir \end{pmatrix}$ as
before,
\[
\z_1((I - tX)x) = ({\bar z}_1 + t(ir{\bar z}_1 + w{\bar z}_2))^n.
\]
Taking the derivative at $t = 0$, we obtain
\[
(X\z_1)(x) = n{\bar z}_1^{n-1}(ir{\bar z}_1 + w{\bar z}_2).
\]
In the same way
\[
(X\z_2)(x) = -n{\bar z}_2^{n-1}({\bar w}{\bar z}_1 + ir{\bar z}_2).
\]
Thus
\[
\<v(x),(Xu)(x)\>_{\cK} \ = \ h^2(x)n(ir2{\bar z}_1^n{\bar z}_2^n \ + \ w{\bar
z}_2^{n+1}{\bar z}_1^{n-1} + {\bar w}{\bar z}_1^{n+1}{\bar
z}_2^{n-1}).
\]
To obtain a lower bound for $L(p)$ we evaluate the absolute value of
this expression at $x_* = (2^{-1/2},2^{-1/2})$ and with $r = 1$, 
$w = 0$.  We obtain $n$.  Thus $L(p) \ge n$.  On the other hand, just
taking absolute values, and using the fact that 
$|z_1|^2 + |z_2|^2 = 1$, we find that for all $x \in G$
\[
|\<v(x),(Xu)(x)\>_{\cK}| \le nh^2(x)(r2|z_1z_2|^n +
|w||z_1z_2|^{n-1}).
\]
The maximum of the right-hand side for $r^2 + |w|^2 = 1$ is
\[
nh^2(x)(4|z_1z_2|^{2n} + |z_1z_2|^{2(n-1)})^{1/2}.
\]
Again using $|z_1|^2 + |z_2|^2 = 1$ we see easily that the maximum
value of $h^2(x)$ is $2^{n-1}$ while the maximum value of $|z_1z_2|$
is $1/2$, and from this we obtain $L(p) \le n\sqrt{2}$.  We summarize
our findings with:

\begin{proposition}  
\label{prop16.4}
For $u$ defined by $u(x) = h(x)(z_1^n,z_2^n)'$ for $n > 0$, or $u(x) =
h(x)({\bar z}_1^n,{\bar z}_2^n)'$ for $n < 0$, and for $p$ the
corresponding projection-valued function, we have
\[
|n| \le L(p) \le |n|\sqrt{2}.
\]
\end{proposition}

It is interesting to notice that as soon as $|n| \ge 2$ the Lipschitz
constant for these projections, coming from standard module frames
with only two elements, is strictly larger, and more rapidly
increasing with $n$, then the projections considered in
Theorem~\ref{th16.2}, which come from standard module frames with
$n+1$ elements.  It is reasonable to guess that these projections have
minimal Lipschitz constant for their two situations, but I have not
examined this matter.

The results of Sections \ref{sec4} and \ref{sec6} can now be applied
in the same way as discussed early in this section for the particular
case of $|n| = 1$.


\section{Appendix A: Path-length spaces}

The purpose of this appendix is to prove:

\begin{proposition}
\label{propA1}
Let $(X, \rho)$ be a compact metric space which is a path-length
space, and let $\cH$ be a real Hilbert space. Let $u:X \to \cH$ be a
continuous function such that $\|u(x)\| = 1$ for all $x \in X$, and
define the projection-valued function $p$ by $p(x) = \<u(x), \
u(x)\>_0$. Then $L(p) = L(u)$.
\end{proposition}

\begin{proof}
Let $v$ and $w \in \cH$ with $\|v\| = 1 = \|w\|$ and 
$\<v, \ w\>_\cH > 0$. Since $\cH$ is over $\bR$, 
\[
\<v, \ w\>_\cH = 1 - (1/2)\|v - w\|^2,
\]
so that from Proposition \ref{prop8.2}
\[
\|\<v, \ v\>_0 - \<w, \ w\>_0\| =
\|v - w\|(1 - (1/4)\|v - w\|^2)^{1/2}.
\]
Thus, since $u$ is uniformly continuous, given any $\e > 0$ we can
find $\d > 0$ small enough that whenever $\rho(x, y) < \d$ then
$\|u(x) - u(y)\| \leq (1+\e)\|p(x) - p(y)\| \leq (1+\e)L(p)\rho(x,y)$.
Thus we need the following lemma, which is essentially known --- see
1.8bis of \cite{Grm2}.

\begin{lemma}
\label{lemA1}
Let $(X, \rho)$ be a path-length metric space, let $B$ be a Banach
space, and let $f$ be a $B$-valued function on $X$. Suppose that there
are constants $K$ and $d > 0$ such that for any $x, y \in X$ with
$\rho(x, y) < d$ we have $\|f(x) - f(y)\| \leq K\rho(x, y)$. Then
$L^\rho(f) \leq K$.
\end{lemma}

\begin{proof}
Let $x, y \in X$ and $\e > 0$ be given, with $\e < d$. Let $\g$ be a
path in $X$ from $x$ to $y$, with domain $I = [0, \ t_*]$, whose
length is 
$\leq ( 1 \ + \ \e/K)\rho(x, y)$. Then $\g$ is uniformly continuous,
so there is a $\d > 0$ such that if $s, t \in I$ with 
$|s \ - \ t| < \d$ then $\rho(\g(s), \ \g(t)) < \e < d$. Choose
$\{t_j\}_{j=0}^m$ in $I$ such that $t_0 = 0$, \ $t_m = t_*$, \
$t_{j+1} > t_j$, and $|t_{j+1} \ - \ t_j| < \d$ for all $j$. Then
\begin{eqnarray*}
\|f(y) \ - \ f(x)\| \ &\leq& \ 
\sum_{j=0}^{m-1} \|f(\g(t_{j+1})) \ - \ f(\g(t_j))\|  \\
&\leq& K\sum \rho(\g(t_{j+1}), \ \g(t_j)) \\ 
&\leq& 
K(1 + \e/K)\rho(x, \ y) \ = \ (K \ + \ \e)\rho(y, \ x).
\end{eqnarray*}
Since $\e$ is arbitrary, we obtain $\|f(x) - f(y)\| \leq K\rho(x, y)$,
as desired. (Actually, we only need that $B$ be a metric space for
the above proof to work.)
\end{proof}
\end{proof}

The above proposition and lemma are false if the path-length
assumption is omitted. An example pertinent to our earlier
constructions of vector bundles can be produced as follows. In $\bR^2$
consider an equilateral triangle with base the unit interval 
$I = [0, \ 1]$ on the $x$-axis of $\bR^2$. Let $X$ be obtained from
the triangle by removing a very small open ball about the vertex
opposite to $I$. Thus $X$ is homeomorphic to a closed interval. Let
$\rho$ be the restriction to $X$ of the Euclidean metric on $\bR^2$,
rather than the evident path-length metric on $X$. Let $u:X \to \bR^2$
be defined by $u(t, 0) = (\cos(\pi t), \ \sin(\pi t))$ on $I$, and by
taking $u$ to be continuous and constant on each of the two ``legs''
of $X$. Thus on one leg $u$ has value $(1, 0)$ while on the other it
has value $(-1, 0)$. Since the two end-points of $X$ are very close
together, $L^\rho(u)$ is very large. But one can find a very small $d$
for which the hypotheses of the lemma are satisfied with a smallish
$K$. On the other hand, when we set $p(x) = \<u(x), u(x)\>_0$, then
$p$ has the same value on the two legs of $X$, and $L(p) < L(u)$.

\providecommand{\bysame}{\leavevmode\hbox to3em{\hrulefill}\thinspace}
\providecommand{\MR}{\relax\ifhmode\unskip\space\fi MR }


\providecommand{\MRhref}[2]{%
  \href{http://www.ams.org/mathscinet-getitem?mr=#1}{#2}
}
\providecommand{\href}[2]{#2}


\end{document}